# Étienne Bézout : Analyse algébrique au siècle des Lumières

## Liliane ALFONSI

**Résumé :** Le but de cet article, à travers l'étude des travaux en analyse algébrique finie d'Étienne Bézout (1730-1783), est de mieux faire connaître ses résultats, tels qu'il les a effectivement trouvés, et de mettre en valeur aussi bien les points de vue novateurs que les méthodes originales, mis en œuvre à cet effet. L'idée de ramener le problème de l'élimination d'une ou plusieurs inconnues à l'étude d'un système d'équations du premier degré, son utilisation inhabituelle des coefficients indéterminés qu'il ne calcule pas mais dont seuls l'existence et le nombre l'intéressent, une façon très personnelle de trouver la résultante de deux équations et enfin l'idée d'étudier, dans leur ensemble, les sommes de produits de polynômes, font partie des approches et des démarches de Bézout que nous nous proposons d'exposer dans cette étude.

**Abstract :** The topic of this paper is, on the one hand to introduce algebraic analysis results of Étienne Bézout (1730- 1783) not as we know them today but as he found them in his time, and on the other hand to emphasize his innovating viewpoints. We will be concerned with Bezout special way of reducing elimination for any degree systems to finding conditions for linear systems solutions, with his typical use of indeterminate coefficients that he doesn't compute but looks only for existence and number, with his idea to work on set of polynomials products sums, and with a very personal method to found two equations resultant.

**Mots clefs :** Analyse algébrique, Bézoutien, coefficients indéterminés, déterminant, élimination, équations algébriques, géométrie algébrique, linéarité, Bézout, Cramer, Euler, Sylvester.

**Key-words:** Algebraic analysis, algebraic equations, algebraic geometry, Bezoutiant, elimination, determinant, indeterminate coefficients, linearity, Bézout, Cramer, Euler, Sylvester.

**Classification mathématique par sujets (AMS 2000) :** 01A50, 12D05, 1403, 15A15

## I. Introduction

Le nom d'Étienne Bézout (1730-1783) est passé à la postérité grâce à deux théorèmes - l'un sur le nombre de points d'intersection des courbes et des surfaces, l'autre, plus connu sous le nom d'« identité », sur les polynômes premiers entre eux - et à son cours de mathématiques, qui a servi de référence à plusieurs générations d'élèves.

Jusque-là, ces résultats et ce cours, sortis de leur contexte, n'avaient suscité aucune étude approfondie sur Étienne Bézout, véritable inconnu célèbre[1]. Sa vie, sa personnalité, aussi bien que son œuvre mathématique et didactique n'avaient fait l'objet que de quelques courts écrits ([Condorcet 1783], [Grabiner 1970]) traitant rapidement et partiellement de l'un ou l'autre de ces aspects, quelquefois de façon erronée ([Vinot 1883], [Lhuillier 1886], [Petit 1930]). C'est à cette carence que nous avons voulu remédier en écrivant une biographie scientifique [Alfonsi 2005] à laquelle on pourra se reporter.

Étienne Bézout est né le 31 mars 1730 à Nemours[2]. Après l'obtention de la maîtrise es arts, il s'installe à Paris vers 1750 comme maître de mathématiques et fréquente rapidement le monde des mathématiciens et des académiciens français [Alfonsi 2005, p. 19-34]. Influencé par d'Alembert, il travaille dans un premier temps (1755-1760) sur des problèmes de mécanique et de calcul intégral[3] et dans ce dernier domaine, on peut déjà remarquer [*Ibid.* p. 38-41, 45-47] son choix d'une approche très algébrique pour résoudre des questions de calcul infinitésimal. Son goût pour la théorie des équations algébriques va se préciser à partir de 1762 et, ses responsabilités très lourdes d'organisateur et d'examinateur d'écoles militaires (Marine à partir de 1764, Artillerie à partir de 1768) l'obligeant par manque de temps à faire un choix, il va se consacrer entièrement à ce domaine.

L'étude de ses travaux nous a amenés à nous interroger sur la spécificité de la démarche mathématique d'Étienne Bézout, académicien, enseignant et chercheur au siècle des Lumières et c'est ce sujet que nous aborderons dans cet article : quels sont exactement les résultats trouvés ? Quelles sont ses problématiques et ses méthodes ? Quelle est sa part d'originalité et d'innovation dans son approche des problèmes et dans les moyens mis en œuvre pour les résoudre ? Dans quel contexte intellectuel et professionnel ses recherches eurent-elles lieu ? Telles sont les questions auxquelles nous allons tâcher d'apporter des éléments de réponse. Nous allons nous intéresser ici, uniquement, aux travaux de Bézout en analyse algébrique finie[4], domaine auquel il a consacré l'essentiel de sa recherche[5] et où son originalité apparaît le mieux.

---

[1] Une illustration de cela est l'hésitation persistante sur l'orthographe même de son nom, Bezout ou Bézout. Nous avons choisi d'orthographier Bézout (avec un accent sur le e et un t) pour respecter son propre choix : cette orthographe est celle de la signature d'Étienne après 1765 et de ses textes imprimés après 1770. Lui-même et les siens signaient Bezout avant 1765, et c'est ainsi qu'est écrit son nom dans les textes imprimés d'avant 1769. Si dans les archives notariales, le nom est toujours Bezout, en revanche dans les registres paroissiaux, on trouve deux écritures : Bezou (le plus souvent) ou Bezout. Cela ne semble pas avoir paru important au propre père d'Étienne, qui signe Bezout des actes dans lesquels son nom est écrit sans t.
[2] Une erreur, sans doute typographique, sur son année de naissance (1739 au lieu de 1730) figure dans l'article du *DSB* [Grabiner 1970]
[3] Voir [Bézout 1760, 1763]

Bézout déplorait d'ailleurs le désintérêt de la majorité des mathématiciens de l'époque pour l'analyse algébrique finie au profit de l'analyse infinitésimale[6] :

« L'application de l'analyse algébrique aux différentes questions qui sont du ressort des Mathématiques, se fait presque uniquement à l'aide des équations. L'analyse infinitésimale également attrayante & importante par les objets nombreux et utiles auxquels on a vu qu'elle pouvait être appliquée, a entraîné tout l'intérêt & tous les efforts; & l'analyse algébrique finie semble, à compter de cette époque, n'avoir été regardée que comme une partie sur laquelle ou il ne restait plus rien à faire, ou dans laquelle ce qui restait à faire, n'aurait été que de vaine spéculation » [Bézout 1779, p. *j-ij*].

Il défendait son choix en faisant référence à plus célèbre que lui :

« Néanmoins la nécessité de perfectionner cette partie [l'analyse algébrique finie], n'a pas échappé à ceux à qui l'analyse infinitésimale est le plus redevable : on a vu que celle-ci même avait besoin que la première fut perfectionnée. Entre plusieurs analystes très distingués, les célèbres MM. Euler & de la Grange ont donné sur cette matière des mémoires qui ne renferment ni moins de profondeur, ni moins de sagacité que les autres productions de ces illustres analystes. » [*Ibid.*]

---

[4] « L'analyse est divisée, par rapport à son objet, en analyse des quantités finies, & analyse des quantités infinies. L'analyse des quantités finies est ce que nous appelons autrement arithmétique spécieuse ou algèbre. L'analyse des quantités infinies […] est celle qui calcule les rapports des quantités qu'on prend pour infinies ou infiniment petites » [Encyclopédie I, 1751, article « analyse »]. L'analyse algébrique finie est donc essentiellement l'étude des équations polynomiales, tandis que l'analyse infinitésimale comprend le calcul différentiel et intégral. Cependant les objets comportant un nombre infini d'opérations, comme les séries, étaient aussi traitées alors de manière algébrique et on a pu parler d'« analyse algébrique » pour qualifier l'ensemble de l'analyse de l'époque, notamment la présentation de Lagrange (voir [Fraser 1987, 1989]).
[5] Pour les autres travaux de Bézout (dynamique, calcul intégral) voir [Alfonsi 2005, p. 35-51]
[6] Les citations sont toujours données avec l'orthographe originale. Nous avons seulement changé parfois, pour une meilleure compréhension, les majuscules et les accents.

Son œuvre algébrique comprend essentiellement cinq écrits : trois mémoires présentés à l'Académie des sciences en 1762, 1764 et 1765, le volume d'Algèbre, datant de 1766, de son cours de mathématiques et son livre *Théorie générale des équations algébriques* publié en 1779. Nous n'aborderons pas ici le mémoire de 1765 qui complète celui de 1762 mais n'apporte rien de nouveau par rapport à notre étude[7]. Nous analyserons ses textes mathématiques et nous les replacerons dans le contexte de travaux antérieurs ou contemporains. D'autre part la prise en compte du milieu et de l'histoire personnelle d'Étienne Bézout, nous permettra d'avancer des hypothèses ou des explications pour certains de ses choix, aussi bien dans les sujets traités que dans la façon dont il les expose et les lieux où il choisit de les publier.

## II.   Le mémoire de 1762 sur la résolution algébrique des équations

Le 20 janvier 1762, il présente à l'Académie des sciences le mémoire « Sur plusieurs classes d'équations de tous les degrés qui admettent une solution algébrique ». Ce travail représente un tournant dans ses recherches qui, à partir de là, ne traiteront plus que d'un seul sujet, les équations : d'abord la résolution algébrique des équations à une seule variable, puis, ce qui est, on le verra, lié, l'élimination de toutes les inconnues sauf une, pour des systèmes d'équations à plusieurs variables. C'est dans ce dernier domaine que pourra s'affirmer l'originalité de sa réflexion.

« Quelque important que soit dans les différentes parties des Mathématiques, la résolution algébrique générale des équations de tous les degrés, nous ne sommes encore guère plus avancés à cet égard, qu'on ne l'était du temps de Descartes : les équations des deuxième,

---

[7] Pour l'étude du mémoire de 1765 [Bézout 1768], voir [Alfonsi 2005, p. 75-91]

troisième et quatrième degrés, sont les seules dans lesquelles on ait pu jusqu'à présent assigner la valeur algébrique générale des racines. » [Bézout 1764a, p. 17]

Voici le constat par lequel Bézout commence son écrit. Il commente longuement l'ouvrage d'Euler de 1732 et fait part de son accord avec une conjecture de celui-ci portant sur l'équation polynomiale générale de degré $n$ à coefficients réels, dont il a annulé de façon standard le terme de degré $n$-1 : « Soit l'équation $x^n + a_{n-2}x^{n-2} + a_{n-3}x^{n-3} + ... + a_1 x + a_0 = 0$, ses racines sont de la forme $x = \sqrt[n]{y_1} + \sqrt[n]{y_2} + ..... + \sqrt[n]{y_{n-1}}$, où les $y_i$ sont les racines d'une équation[8] de degré $n$-1 » [Euler 1738, p. 8].

Pour résoudre une équation à une inconnue de degré $n$ quelconque, l'idée de Bézout dans ce mémoire est la suivante : « On peut […] prendre arbitrairement une équation à deux termes en y, & du degré de la proposée, comparer à cette équation une équation en $x$ & $y$, telle que de cette comparaison il résulte une équation du même degré que la proposée, avec laquelle on la comparera terme à terme. » [Bézout 1764a, p. 22]

Il va donc considérer l'équation générale $(A): x^n + mx^{n-1} + px^{n-2} + qx^{n-3} + rx^{n-4} + ... + M = 0$, de degré $n$ quelconque, comme le résultat du système $\begin{cases} (B): y^n + h = 0, \\ (C): y = \dfrac{x+a}{x+b} \end{cases}$.

Ceci lui servira de fil conducteur pour traiter les trois problèmes qui structurent cet écrit :

- résoudre l'équation générale de degré 3 en la réduisant à une équation du même degré à deux termes ;

- trouver les conditions qui réduiraient une équation de degré quelconque à une équation du même degré à deux termes ;

- trouver des équations résolubles par la somme de 2, 3, 4, etc., « radicaux du degré de ces équations » [*Ibid*., p. 33].

---

[8] Équation sans second terme, c'est à dire sans terme de degré n-2, par annulation standard du second terme.

En appliquant sa méthode au degré 3, Bézout traite le premier point, à quelques cas particuliers près [Alfonsi 2005, p. 70]. Mais, depuis les algébristes italiens du XVIᵉ siècle, on savait résoudre les équations du troisième degré, y compris par des méthodes très proches de la sienne comme celle de Tschirnaüs exposée en 1683 dans les *Acta Eruditorum* de Leipzig – méthode que Bézout ne connaissait pas, semble-t-il, en 1762 [*Ibid.*, p. 67].

Nous allons donc nous intéresser surtout aux deux autres problèmes. Pour trouver les conditions qui réduiraient une équation de degré quelconque à une équation du même degré à deux termes, les notations étant celles données plus haut, Bézout élimine $y$ entre $(B)$ et $(C)$, et obtient l'équation [9] $(D) : \sum_{k=0}^{n} C_n^k x^{n-k} \dfrac{a^k + h b^k}{1+h} = 0$.

Il considère ensuite, ce qui est toujours possible à un changement de variable près, que $m=0$, et donc par identification, que $h = -a/b$, ce qui lui permet d'écrire que $(D)$ donne l'équation

$(E) : x^n - C_n^2 x^{n-2} ab - C_n^3 x^{n-3} ab(a+b) - C_n^4 x^{n-4} ab(a^2 + ab + b^2) - .. - ab(a^{n-2} + a^{n-3}b + .. + b^{n-2}) = 0$

En identifiant les coefficients des termes en $x^{n-2}$ et $x^{n-3}$ de $(A)$ et de $(E)$, c'est à dire,

$\begin{cases} C_n^2 ab = -p \\ C_n^3 ab(a+b) = -q \end{cases}$, il obtient $a$ et $b$ comme solutions de $\quad X^2 - \dfrac{3q}{(n-2)p} X - \dfrac{p}{C_n^2} = 0$ .

Il n'envisage que le cas général, ne s'occupant pas des exceptions [*Ibid.*, p. 71]. Ceci lui permet de trouver « les conditions qui réduiraient une équation de degré quelconque à une équation du même degré à deux termes » [Bézout 1764a, p. 27] par identification des coefficients restants de *(A)* et *(E),* dans laquelle il remplace $a$ et $b$ par leurs valeurs en fonction de $p$ et $q$.

Quand l'équation générale de degré $n$ vérifie les conditions qui la réduisent à $y^n + h = 0$, il trouve, en remplaçant $h$ par $-a/b$ et $y$ par la valeur qui en découle, que les

---

[9] Bézout n'utilise pas la notation $C_n^k$, mais l'écriture développée en quotient de produits d'entiers. C'est pour faciliter l'écriture et la lecture que nous utilisons des notations plus actuelles.

racines de l'équation de degré $n$ s'écrivent comme la somme de $n$-1 racines $n$-ièmes :

$$x = \sqrt[n]{a^{n-1}b} + \sqrt[n]{a^{n-2}b^2} + \ldots + \sqrt[n]{ab^{n-1}}$$ , $a$ et $b$ étant les solutions d'une équation du $2^{nd}$ degré.

« On voit par-là, écrit-il, que, dans une équation de degré quelconque et sans second terme, les coefficiens du troisième et du quatrième terme étant tels qu'on voudra, si les coefficiens des autres termes sont tels qu'ils résultent de la comparaison des deux équations $(E)$ et $(F)$ [$(F)$ étant l'équation générale de degré $n$, sans second terme], cette équation sera résoluble et aura pour racines la somme de $n$-1 moyennes proportionnelles entre les deux racines d'une équation du deuxième degré. » [Bézout 1764a, p. 27]

Bézout retrouve ainsi, l'écriture conjecturée par Euler en 1732, mais avec des conditions très contraignantes sur les coefficients de l'équation de départ ; il ne parvient pas, bien sûr, à la résolution du cas général de degré $n$, mais il ajoute une classe d'équations résolubles, à celles déjà trouvées par Moivre et Euler.

Dans la troisième partie, Bézout revient à la démarche inverse de la précédente, qui fut aussi celle d'Euler dans son mémoire de 1732 : en partant d'une somme de racines $n$-ièmes, chercher les équations dont cette somme peut être racine.

Posant $x = \sqrt[n]{a^{n-1}b} + \sqrt[n]{a^{n-2}b^2}$ et après avoir traité les cas particuliers, $n = 3, 4, 5, 6, 7, 8$, il parvient par induction, pour tout $n$ à l'équation :

$$x^n = a^{n-1}b \pm a^{n-2}b^2 + na^{n-2}bx + n\frac{n-3}{2}a^{n-3}bx^2 + n\frac{n-4}{2}\frac{n-5}{2}a^{n-4}bx^3 + \ldots$$

$$\ldots + n\frac{n-p}{2}..\frac{n-2p+3}{2}a^{n-p}bx^{p-1} + \ldots$$

« suite, dit-il, qu'on doit pousser seulement jusqu'au terme dont le coefficient devient zéro[10] , & dans laquelle le signe $+$ du terme $a^{n-2}b^2$ est pour les degrés impairs, & le signe $-$ pour les degrés pairs ; donc toute équation qui pourra se rapporter à cette dernière, sera résoluble en faisant $x = \sqrt[n]{a^{n-1}b} + \sqrt[n]{a^{n-2}b^2}$ » [*Ibid.* p. 37].

---

[10] Il ne fait donc pas de calcul sur des séries formelles.

Il traite ensuite le cas $x = \sqrt[n]{a^{n-2}b^2} + \sqrt[n]{a^{n-3}b^3}$ , pour $n$ allant de 5 à 9, mais ne parvient pas à induire une équation pour tout $n$, car « les coefficients suivent une loi moins uniforme ». De même, il travaille sur d'autres cas particuliers de $x = \sqrt[n]{a^{n-m}b^m} + \sqrt[n]{a^{n-m-1}b^{m+1}}$ jusqu'à $n=7$.

À la lecture de ce mémoire, on voit apparaître les éléments qui vont déterminer l'orientation future du travail de recherche d'Étienne Bézout. Tout d'abord la résolution des équations à une variable par l'élimination d'une inconnue $y$ dans un système de deux équations à deux inconnues $x, y$, pose le problème du degré de l'équation « réduite »[11].

D'autre part dans ce travail, Bézout utilise beaucoup la méthode des coefficients indéterminés, méthode classique que d'Alembert attribue à Descartes : « La méthode des *coefficiens* indéterminés est une des plus importantes découvertes que l'on doive à Descartes » [*Encyclopédie,* art. « coefficient », t. 3, 1753]. On la trouve ainsi décrite dans *La Géométrie* de ce dernier : « Mais je veux bien, en passant, vous avertir que l'invention de supposer deux équations de mesme forme, pour comparer séparément tous les termes de l'une à ceux de l'autre, & ainsi en faire naistre plusieurs d'une seule, dont vous avez vu icy un exemple, peut servir à une infinité d'autres problesmes & n'est pas l'une des moindres de la méthode dont je me sers » [Descartes 1637, p. 351].

D'Alembert lui-même l'expose ainsi :

« Cette méthode [...]consiste à supposer l'inconnue égale à une quantité dans laquelle il entre des *coëfficiens* qu'on suppose connus, & qu'on désigne par des lettres ; on substitue cette valeur de l'inconnue dans l'équation ; & mettant les uns sous les autres les termes homogènes,

---

[11] À l'époque, l'équation obtenue par élimination d'inconnues dans un système s'appelle plutôt la « réduite », si le système est celui que l'on écrit en vue de résoudre une équation à une inconnue (la réduite d'une équation), et plutôt la « résultante », si le système est quelconque (la résultante d'un système). Mais dans les textes, ces deux expressions sont parfois utilisées indifféremment dans ces deux situations.

on fait chaque *coefficient* = 0, & **on détermine par ce moyen les *coefficiens* indéterminés**. »
[*Encyclopédie,* art. « coefficient », t. 3, 1753] (souligné par nous)

Il donne un exemple dans le cas d'une équation différentielle dont il cherche la solution sous la forme d'un polynôme du second degré, écrit a priori avec des coefficients indéterminés :

soit l'équation $dy + (by + ax^2 + cx + f)dx = 0$, il suppose $y = A + Bx + Cx^2$, où *A, B*, et *C* sont des coefficients indéterminés, alors $dy = Bdx + 2Cxdx$, et $bydx = bAdx + bBxdx + bCx^2dx$ ; en identifiant les deux formes, d'Alembert obtient le système de trois équations à trois inconnues, $B + bA + f = 0$ , $2C + bB + c = 0$ , $bC + a = 0$ , système qui lui permet de calculer *A, B*, et *C*.

Il est important de rappeler ces définitions car, si en 1762 Étienne Bézout utilise la méthode de façon classique (égaler une quantité connue à une expression contenant des coefficients indéterminés, calculer les coefficients) nous verrons que dans ses travaux ultérieurs, il n'utilisera plus la méthode des coefficients indéterminés de cette façon.

### III.     Le mémoire de 1764 sur l'élimination des inconnues

En 1764, Étienne Bézout présente à l'Académie des Sciences un écrit qui deviendra l'un de ses deux plus importants travaux de recherche (l'autre étant la *Théorie générale des équations algébriques* de 1779), le mémoire « Recherches sur le degré des équations résultantes de l'évanouissement des inconnues et sur les moyens qu'il convient d'employer pour trouver ces équations ». Il le lit au cours des séances des 1[er], 15, 24 et 29 février et cet ouvrage est publié en 1767 dans les *Mémoires de l'Académie royale des sciences* pour 1764.

### 1.     Le problème de l'élimination

Bien que nous ayons déjà parlé d'élimination dans le contexte de la résolution algébrique des équations, nous allons revenir sur cette notion dans un cadre plus général.

Le problème de l'élimination peut s'énoncer ainsi : un certain nombre d'inconnues et de relations polynomiales entre ces inconnues étant donné, y-a-t-il des valeurs de ces inconnues qui vérifient toutes ces relations ? Et dans le cas d'une réponse positive, quelles sont-elles ?

Pour résoudre ce double problème, on cherche à déduire, à partir des relations polynomiales données, une équation dans laquelle ne subsiste plus, au maximum, qu'une seule inconnue. On dit alors que l'on a « éliminé » les autres. Cette équation donnera une condition d'existence, d'où l'on pourra tirer éventuellement les valeurs de l'inconnue restante qui, une fois calculées, permettront de trouver celles des autres inconnues.

Exemple[12] : soient $P(x,y)$ et $Q(x,y)$ deux polynômes en $x$ et $y$ tels que, une fois ordonnés par rapport aux puissances de $x$, on ait le système :
$$\begin{cases} P(x, y) = A(y)x^2 + B(y)x + C(y) = 0 \\ Q(x, y) = A'(y)x^2 + B'(y)x + C'(y) = 0 \end{cases}$$

Éliminer $x$ dans ce système, donne la condition d'existence de solutions en x,

$[A(y)B'(y) - B(y)A'(y)] \, [C'(y)B(y) - B'(y)C(y)] = [A'(y)C(y) - A(y)C'(y)]^2$.

Cette relation est une équation en $y$ - la résultante du système -, dont la résolution permet d'avoir les $x$ correspondants.

Pour mieux situer le mémoire de 1764 d'Étienne Bézout et montrer ce qu'il apporte de nouveau et d'original, nous allons évoquer rapidement les travaux de ses prédécesseurs les plus importants, Newton[13], Euler et Cramer qu'il cite lui même dans son introduction. Comme dans son premier mémoire sur la résolution algébrique des équations, Bézout, de façon très pédagogique, commence par un rapide historique du problème de l'élimination. Il veut, comme il l'écrit, « ne point envelopper dans [son] travail ce qui peut appartenir à d'autres » [Bézout 1767c, p. 289]. L'élimination des inconnues s'était essentiellement portée jusque-là, sur le cas de deux équations à deux inconnues, cas qui correspondait à la recherche des points d'intersection de deux courbes planes, et que l'on résolvait surtout par substitution, ce qui entraînait souvent des facteurs superflus dans l'équation finale.

---

[12] Cet exemple, sous une forme légèrement différente, figure dans le mémoire de Bézout [Bézout 1767c, p. 319]. Il y est traité par la méthode dite plus tard du « Bézoutien », voir *infra*.
[13] Newton avait, dans son *Arithmetica universalis* [Newton 1707] donné directement les formules de la résultante, sans facteurs superflus, pour deux équations de degré $n$ chacune, avec $n=1, 2, 3$ ou 4. Par ailleurs, le nombre de points d'intersection de deux courbes planes était pour lui, de toute évidence et sans démonstration, le produit des degrés des deux courbes. On le voit, dans un de ses écrits rédigé vers 1665 et non publié de son temps [Newton *MP*, Vol.1], pour lui c'est plus un principe qu'un problème. Il affirme qu'« il suffit » d'éliminer une inconnue, et on trouve « autant d'intersections que le rectangle des dimensions des courbes » [*Ibid*, p. 498].

Euler est le premier qui ait donné une démonstration du fait qu'une courbe de degré $m$ et une courbe de degré $n$, se coupaient au plus en $m.n$ points, même si cette preuve contient des points à éclaircir. Il présente en 1748 à l'Académie des Sciences de Berlin, un écrit intitulé « Démonstration sur le nombre de points où deux lignes des ordres quelconques peuvent se couper », qui sera publié en 1750 dans les *Mémoires de l'Académie* pour 1748. Dans son introduction il situe l'état du problème :

« Dans la pièce précédente [Euler 1750a], j'ai rapporté sans démonstration cette proposition, *que deux lignes courbes algébriques, dont l'une est de l'ordre[14] m et l'autre de l'ordre n, se peuvent couper en m.n points*. La vérité de cette proposition est reconnue de tous les géomètres, quoiqu'on doive avouer, qu'on n'en trouve nulle part une démonstration assés rigoureuse. Il y a des vérités générales que notre esprit est prêt d'embrasser aussitôt qu'il en reconnoit la justesse dans quelques cas particuliers : et c'est parmi cette espèce de vérités qu'on peut ranger à bon droit la proposition, dont je viens de faire mention, puisqu'on la trouve vraie non seulement dans quelques ou plusieurs cas, mais aussi dans une infinité de cas différens. Cependant on conviendra aisément que toutes ces preuves infinies ne sont pas capables de mettre cette proposition à l'abri de toutes les objections qu'un adversaire peut former, et qu'il faut absolument une démonstration rigoureuse, pour le réduire au silence. » [Euler 1750b, p. 46]

Il explique ensuite ce qu'il faut entendre par points d'intersection de deux courbes et avance une conception tout à fait moderne, qui se rapproche de l'énoncé actuel du « théorème de Bézout » : « le sens de notre proposition est que le nombre des intersections ne peut jamais être plus grand que *m.n*, quoiqu'il soit souvent plus petit ; et alors on juge, ou que quelques intersections s'éloignent à l'infini, ou qu'elles deviennent imaginaires. De sorte qu'en contant les intersections à l'infini et les imaginaires aussi bien que les réelles, on puisse dire que le nombre des intersections est toujours = *m.n* » [*Ibid.*]. Par ailleurs il écarte le cas où les deux courbes ont une infinité de points communs, donc une composante commune.

---

[14] L'ordre est le degré du polynôme définissant la courbe

Après quelques exemples par lesquels il montre qu'on peut facilement obtenir, dans la résultante, des facteurs qui donnent de fausses racines il en vient à la méthode qu'il va suivre : En partant de deux courbes planes, donc de deux équations à deux inconnues (*x* et *y*) de degrés respectifs *m* et *n*, il choisit d'ordonner chaque équation par rapport aux puissances de *y*. Il considère donc le système suivant :

$$\begin{cases} f(y) = y^m - Py^{m-1} + Qy^{m-2} - Ry^{m-3} + Sy^{m-4} - etc. = 0 \\ g(y) = y^n - py^{n-1} + qy^{n-2} - ry^{n-3} + sy^{n-4} - etc. = 0 \end{cases}$$

où *P, Q, R, S,* etc. *s*ont des polynômes en *x*, de degrés respectifs, 1, 2, 3, etc. jusqu'à *m*, et *p, q, r, s,* etc. *s*ont des polynômes en *x*, de degrés respectifs, 1, 2, 3, etc. jusqu'à *n*.

Si, *A, B, C, D,* etc., sont les *m* racines de la première équation, *f(y)*=0, et *a, b, c, d,* etc., les *n* racines de la seconde, *g(y)*=0, il peut écrire :

$$\begin{cases} f(y) = (y-A)(y-B)(y-C)(y-D)etc. = 0 \\ g(y) = (y-a)(y-b)(y-c)(y-d)etc. = 0 \end{cases}.$$

Alors pour que les deux équations aient une racine commune, il faut et il suffit que l'une des racines *A, B, C, D,* etc., de la première équation soit égale à l'une des racines *a, b, c,* etc., de la seconde, donc il faut et il suffit que

[(A-a)(A-b)(A-c) *etc.*] . [(B-a)(B-b)(B-c).*etc.*] . [(C-a)(C-b)(C-c) *etc.*] . etc. = 0.

Autrement dit,         *g(A).g(B).g(C).g(D)*.etc. = 0 ,

ou,         *f(a).f(b).f(c).f(d)*.etc. = 0

Ces dernières équations, qui sont équivalentes, représentent, pour lui, la résultante cherchée. Comme elles représentent une condition nécessaire et suffisante pour que les deux équations de départ aient une racine commune, il est sûr que cette résultante ne contient pas de facteurs superflus.

Euler est conscient que « les expressions des racines *A, B, C, D*, etc., *a, b, c, d*, etc., sont pour la plupart fort irrationnelles et souvent telles, qu'on ne les peut pas assigner » [*Ibid.* p. 55], et qu'il faut montrer que sa résultante est bien un polynôme en *x*. Il résout le problème en se servant des relations entre les coefficients et les racines d'une équation. Il affirme que l'« on voit aisément », que dans *g(A).g(B).g(C).g(D).*etc. = *0*, les termes obtenus par les produits des différentes puissances des *A, B, C, D*, etc., peuvent être remplacés par des produits de différentes puissances des *P, Q, R*, etc., qui eux sont des polynômes en *x*, mais il ne le montre vraiment, que sur des cas particuliers de degrés 2 et 3.

Il a raison sur ce dernier point – bien qu'il ne l'ait pas montré - car, si on considère sa résultante sous la forme *g(A).g(B).g(C).g(D).*etc. = *0*, il s'agit bien d'un polynôme symétrique en *A, B, C*, etc., donc aussi, nous le savons aujourd'hui, un polynôme en *P, Q, R*, etc., donc un polynôme en *x*.

Par ailleurs, on a $\begin{cases} A + B + C + D + etc. = P \\ a + b + c + d + etc. = p \end{cases}$ avec *P* et *p* chacun de degré 1 en *x*.

De même la somme des produits deux à deux des racines *A, B, C, D*, etc., est de degré 2, de même pour *a, b, c, d*, etc., et plus généralement pour un nombre *k,* compris entre 1 et *m* pour *f*, et entre 1 et *n* pour *g*, la somme des produits *k* à *k* des racines *A, B, C, D*, etc., est de degré *k* en *x,* de même pour *a, b, c, d,* etc. Il conclut de là que les racines *A, B, C, D*, etc., et *a, b, c, d,* etc., peuvent être considérées chacune comme un polynôme de degré 1 en *x*. Étant donné la forme de la résultante trouvée, produit de *m.n* facteurs de degré au plus 1, elle ne peut être que de degré au plus *m.n* en *x*.

On voit tout de suite quel est le point faible majeur de ce raisonnement :

Euler, comme d'Alembert, avait donné en 1746, une démonstration, lacunaire mais sérieuse, du théorème fondamental de l'algèbre, pour un polynôme à une variable à coefficients réels. Il énonce que, tout polynôme à une variable, de degré *n* et à coefficients réels, a *n* racines réelles ou imaginaires.

Dans ce mémoire de 1748, Euler étend le résultat par analogie à un polynôme dont les coefficients sont eux-mêmes polynomiaux, mais aucun résultat, à son époque, ne lui permettait de le faire. En langage actuel, les racines *A, B, C, D,* etc., et *a, b, c, d,* etc., qu'il considère, existent, mais ce sont des éléments de la clôture algébrique de ***R****(X),* c'est à dire de la clôture du corps des fractions rationnelles sur ***R***.

Il est alors difficile de dire, comme il le fait, que : « on pourra regarder chaque racine comme une fonction d'une dimension de *x* ». En effet, d'une part il s'est placé d'emblée dans un cas particulier, car *P, Q, R,* etc., et *p, q, r,* etc., n'ont pas obligatoirement les degrés (1 pour *P* et *p,* 2 pour *Q* et *q,* etc.) qu'il leur impose (les degrés peuvent être inférieurs à son choix), d'autre part, même dans la situation choisie, un simple exemple du second degré montre bien que la somme des racines étant du premier degré et le produit du second, les racines ne sont pas obligatoirement du premier degré, ni mêmes des polynômes[15].

Pourtant, en admettant sa factorisation de *f* et de *g,* sa résultante est bien de degré inférieur ou égal à *m.n* en *x* : cette résultante est un polynôme symétrique en *A, B, C, D,* etc., comme nous l'avons vu, et ce polynôme est d'ordre *n* et de degré *m.n* par rapport à ces variables ; cette même résultante est donc un polynôme de degré *n* par rapport aux coefficients *P, Q, R,* etc. ; donc, sa résultante est de degré inférieur ou égal à *m.n* en *x*.

En conclusion, si la démonstration d'Euler était lacunaire dans l'état des connaissances de l'époque (et il en avait conscience puisqu'il écrivait : « S'il y a dans cette démonstration encore quelque obscurité, cela vient de sa grande généralité » [Euler 1750b, p. 57]), elle peut être établie rigoureusement aujourd'hui.

---

[15] Si l'on prend y²-(x-1)y-(x²+2)=0 comme un polynôme en y, les racines $\frac{1}{2}(x-1\pm\sqrt{5x^2-2x+9})$ ne sont pas des polynômes, et pourtant leur somme est un polynôme de degré 1 et leur produit, un polynôme de degré 2.

En ce qui concerne les travaux d'Euler sur l'élimination, il faut citer aussi le chapitre XIX du tome 2 de l'ouvrage *Introductio in analysin infinitorum* [Euler 1748], ainsi que le mémoire intitulé « Nouvelle méthode d'éliminer les quantités inconnues des équations » [Euler 1766] qu'il présente la même année que Bézout, mais lui, à l'Académie des sciences de Berlin et qui sera publié en 1766. Dans ces deux écrits Euler présente des techniques différentes pour obtenir la résultante mais il ne revient pas sur le calcul de son degré.

Cramer, cité lui aussi par Bézout, a publié à Genève en 1750 son *Introduction à l'analyse des lignes courbes algébriques,* dans lequel, bien qu'il ne semble pas avoir eu connaissance du mémoire d'Euler de 1748, publié lui aussi en 1750, il emploie la même idée (présentée sous une forme différente), pour arriver au même résultat sur les intersections de courbes planes. Comme il le dit dans sa préface, en exposant les différents sujets qu'il va aborder dans son livre, Cramer pense être le premier à l'avoir démontré : « On y voit [...] le nombre des points dans lesquels une ligne d'un ordre donné peut rencontrer une ligne du même ordre, ou d'un autre ordre aussi donné. La règle qui détermine ce nombre est très importante dans la théorie des courbes, plusieurs grands géomètres l'ont supposée, mais personne, que je sache, n'en a donné la démonstration » [Cramer 1750, p. XIII].

Il note[16] le système des deux équations *A* et *B* à deux inconnues :

$$\begin{cases} A.....................x^n - [1]x^{n-1} + [1^2]x^{n-2} - [1^3]x^{n-3} + \&c.......[1^n] = 0 \\ B............(0)x^0 + (1)x^1 + (2)x^2 + (3)x^3 + \&c..............(m)x^m = 0 \end{cases}$$

dans lequel, contrairement à Euler, il se place dans le cas général, considérant que (*m*) est une constante, [1] et (*m-1*) des polynômes en *y* de degré au plus 1, [1²] et (*m-2*) des polynômes en *y* de degré au plus deux, etc., [1ⁿ] est de degré au plus *n* et (0) de degré au plus *m*.

---

[16] L'écriture du système est exactement celle de Cramer (points de suspension y compris). Cramer a l'idée originale d'introduire une numérotation dans les coefficients qui en indique l'emplacement. Ses notations se rapprochent de notations plus modernes et facilitent la compréhension des calculs. Voir Cajori [1928 / 1993, p. 398].

Il part lui-aussi de la factorisation linéaire d'une équation, mais ne l'utilise que pour l'équation *A*, notant *a, b, c, d,* etc., ses *n* racines. Si nous notons par $f(x)=0$ l'équation *A*, et par $g(x)=0$ l'équation *B,* nous pouvons noter alors, comme Euler, l'équation résultante qu'il obtient :

$$g(a).g(b).g(c).g(d).\text{etc.} = 0.$$

Il ne remarque pas, contrairement à Euler, que cette équation ne contient sûrement aucun facteur superflu. À partir de là, sa méthode diffère. Au lieu de garder la résultante sous l'une des formes que nous avons vues chez Euler, il la développe, autrement dit, il effectue tous les produits $g(a).g(b).g(c).g(d).$etc.

En revanche, il ne fait pas de remarque sur la nature des racines *a, b, c, d,* etc., mais déplore simplement que ces racines soient inconnues « lorsque l'équation A est d'un degré trop élevé pour que l'Algèbre en puisse donner la solution » [*Ibid*. p.662]. C'est pour cela qu'il utilise les coefficients de A. En partant des relations entre les coefficients et les racines, et grâce à des jeux d'écriture facilités par ses notations qui lui permettent d'évaluer facilement les exposants de *y*, il montre que le degré de la résultante en *y* ne peut pas dépasser *m.n.*

La méthode de Cramer part donc du même principe que celle d'Euler, mais, si elle est moins élégante - beaucoup plus longue, laborieuse, calculatoire et basée sur des astuces de notations -, elle est plus satisfaisante pour le calcul du degré de la résultante, car Cramer se place dans le cas général et démontre vraiment que le degré ne peut dépasser *m.n.*

**2.      Le calcul du degré de la résultante de deux équations à deux inconnues**

Revenons maintenant au travail d'Étienne Bézout. Dès l'introduction celui-ci fait le lien avec son premier travail sur la résolution algébrique des équations, écrit deux ans plus tôt : « Les recherches [il s'agit de travaux sur l'élimination] dont je vais exposer les résultats dans ce mémoire, doivent naissance à *celles dont je continue de m'occuper sur la résolution* *algébrique des équations* [souligné par nous][17] » [Bézout, 1767c, p. 92]. Cependant, Bézout a découvert l'intérêt propre et plus général de la théorie de l'élimination, puisqu'il ajoute :

« Si les méthodes d'élimination n'avoient d'autre utilité que leur application à la résolution algébrique des équations, je me serois contenté de ce qui peut avoir rapport à ce dernier objet & je l'aurois réuni avec ce que j'ai pu trouver jusqu'à présent sur cette matière ; mais ces méthodes ont une application beaucoup plus étendue & telle qu'elles deviennent indispensables dans tous les problèmes où il y a plus d'une inconnue. En effet, si on a des méthodes pour résoudre, par approximation, les problèmes déterminés lorsqu'on n'a qu'une seule équation, on n'en a pas de même pour les résoudre par approximation lorsque les relations des inconnues, qui en font l'objet, restent, pour ainsi dire, dispersées dans plusieurs équations ; ainsi quand même on seroit condamné pour toujours à résoudre par approximation, les méthodes d'élimination n'en seroient pas moins indispensables » [*Ibid.*].

On s'aperçoit dans ce dernier passage, qu'il envisage l'irrésolubilité des équations de degré 5 et au-delà. En tout cas, amené au problème de l'élimination par la résolution des équations, il la considère maintenant comme un sujet à part entière et cela va conduire ses recherches vers ce domaine de l'algèbre dans lequel il excellera.

Bézout expose une idée nouvelle et personnelle pour traiter ce problème : « Je réduis, dans ce mémoire, tout le travail de l'élimination, à quelques degrés que montent les équations, je le réduis, dis-je, à éliminer des inconnues au premier degré. » [Bézout 1767c, p. 291]

Cette méthode, caractéristique de Bézout, consiste en ramenant – nous verrons de quelle façon - la résolution d'un système quelconque à celle d'un système linéaire homogène, a écrire l'équation résultante comme condition d'existence de solutions non nulles à ce système linéaire, c'est-à-dire à annuler son déterminant[18]. On peut dès lors comprendre la progression de son ouvrage qui se développe en plusieurs étapes :

- améliorer la formation du déterminant d'un système linéaire,

- ramener le calcul de la résultante à celui du déterminant d'un système linéaire homogène ;

- appliquer cette méthode, d'abord, à la recherche du degré de la résultante de 2 équations à 2 inconnues, puis de $n$ équations à $n$ inconnues, avec $n>2$ ;

- puis l'appliquer pour mettre au point une méthode performante de calcul de la résultante, dans le cas de deux équations à deux inconnues (celle qui a conduit à l'introduction du Bézoutien).

---

[17] Cette phrase, écrite en 1764, explique son mémoire de 1765 sur la résolution des équations, sujet qu'il n'avait pas abandonné.

Une des grandes originalités du travail d'Étienne Bézout sur l'élimination est là : pour lui, la résultante est toujours donnée par le déterminant d'un système linéaire. Cette idée est tout à fait nouvelle et sera son fil conducteur.

Pour la formation du déterminant, tout en rendant hommage à Cramer, Bézout pense avoir apporté quelques améliorations à l'application de ses formules :

« Il n'y a encore que fort peu de temps qu'on a une méthode pour trouver la valeur des inconnues dans les équations du $1^{er}$ degré d'une manière simple & sans que cette valeur soit compliquée de quelque facteur inutile; [...] M. Cramer a donné une règle générale pour les exprimer toutes débarrassées de ce facteur : j'aurais pu m'en tenir à cette règle[19]; mais l'usage

---

[18] Bézout n'emploie pas le mot déterminant, ni la notation que nous employons aujourd'hui (voir [Muir 1906] et [Knobloch 1994]). Quand cela s'avèrera nécessaire par commodité et pour une compréhension actuelle des calculs et des résultats, nous les emploierons tout de même par la suite, en le signalant chaque fois.

[19] La règle de Cramer qui correspond à la formation de ce que nous appelons aujourd'hui le déterminant est la suivante : « Soient plusieurs inconnues z, y, x, &c., et autant d'équations

$A^1 = Z^1 z + Y^1 y + X^1 x + V^1 v + \&c.$

$A^2 = Z^2 z + Y^2 y + X^2 x + V^2 v + \&c.$

$A^3 = Z^3 z + Y^3 y + X^3 x + V^3 v + \&c.$ où les lettres $A^1, A^2, A^3, A^4, \&c.$, ne marquent pas comme à

$A^4 = Z^4 z + Y^4 y + X^4 x + V^4 v + \&c.$

$\&c.$

l'ordinaire les puissances de $A$, mais le premier membre supposé connu de la première, seconde, troisième, quatrième, &c. équation. [...] Le nombre des équations et des inconnues étant $n$, on trouvera la valeur de chaque inconnue en formant $n$ fractions dont le dénominateur commun a autant de termes qu'il y a de divers m'a fait connaître que quoiqu'elle soit assez simple, quant aux lettres, elle ne l'est pas de même à l'égard des signes lorsqu'on a au-delà d'un certain nombre d'inconnues à calculer ; j'ai donc cru devoir revenir sur cet objet. » [Bézout 1767c, p. 291]

Étienne Bézout va donc reprendre, à sa façon, la formation du déterminant d'un système linéaire et cela aura une grande importance dans son œuvre. Au lieu de compter pour chaque terme les « dérangements » et en déduire son signe - comme le propose Cramer - il indique qu'il « réduit le travail à n'exiger d'autre attention que celle qu'il faut pour écrire des lettres ». Comme il le dit au début de son lemme 1 :

« Si l'on a un nombre *n* d'équations du premier degré qui renferment chacune un pareil nombre d'inconnues, *sans aucun terme absolument connu* [souligné par nous, il n'y a donc pas de terme constant à gauche du signe égal, et seulement zéro à droite de ce signe, comme on le constate dans les exemples que Bézout donne ensuite], on trouvera par la règle suivante la relation que doivent avoir les coefficients de ces inconnues *pour que toutes ces équations aient lieu* [souligné par nous, Bézout veut donc des solutions non nulles] » [Bézout 1767c, p. 292].

Le problème qui intéresse Bézout, n'est pas, contrairement à Cramer, de résoudre un système quelconque, mais, de trouver la condition pour qu'un système homogène ait des solutions non toutes nulles. Il énonce alors, *sans aucune démonstration*, sa règle dont on peut supposer qu'il l'a obtenue par induction[20] :

« Soient *a, b, c, d*, &c., les coëfficiens de ces inconnues dans la première équation.

*a', b', c', d'*, &c., les coëfficiens de ces inconnues dans la seconde équation

---

arrangements de *n* choses différentes. Chaque terme est composé des lettres ZYXV, &c., toujours écrites dans le même ordre, mais auxquelles on distribue, comme exposants, les *n* premiers chiffres rangés en toutes les manières possibles. […] On donne à ces termes les signes + ou -, selon la règle suivante. Quand un exposant est suivi dans le même terme, médiatement ou immédiatement, d'un exposant plus petit que lui, j'appellerai cela un *dérangement*. Qu'on compte, pour chaque terme, le nombre des dérangements : s'il est pair ou nul, le terme aura le signe + ; s'il est impair, le terme aura le signe -. Le dénominateur étant ainsi formé, […] » [Cramer 1750, p. 657-659]

[20] On peut penser qu'il s'est inspiré de la règle de Cramer, donnée dans la note ci dessus, règle elle-même obtenue par induction à partir des cas *n* = 1, 2, et 3 [Cramer 1750, p. 657-658]

*a″, b″, c″, d″*, &c., ceux de la troisième & ainsi de suite.

[…] Formez les deux permutations *ab* & *ba* & écrivez *ab-ba* ; avec ces deux permutations & la lettre *c* formez toutes les permutations possibles en observant de changer de signe toutes les fois que *c* changera de place dans *ab* & la même chose à l'égard de *ba*; vous aurez

*abc-acb+cab-bac+bca-cba.*

Avec ces six permutations et la lettre *d* formez toutes les permutations possibles, en observant de changer de signe à chaque fois que *d* changera de place dans un même terme [...].et ainsi de suite jusqu'à ce que vous ayez épuisé tous les coëfficiens de la première équation. Alors conservez les lettres qui occupent la première place; donnez à celles qui occupent la seconde, la même marque qu'elles ont dans la seconde équation; à celles qui occupent la troisième, la même marque qu'elles ont dans la troisième équation, & ainsi de suite; égalez enfin le tout à zéro et vous aurez l'équation de condition cherchée » [*Ibid.*]

En termes actuels, Bézout donne un algorithme pour écrire un déterminant d'ordre quelconque et il affirme qu'un système homogène a des solutions non nulles si et seulement si son déterminant est égal à zéro. Il met ensuite ces conditions sous la forme :

*ab'-a'b*=0,

*(ab'-a'b)c″+(a″b-ab″)c'+(a'b″-a″b')c*=0     &c…

« Cette nouvelle forme a deux avantages, affirme-t-il, le premier, de rendre les substitutions à venir, plus commodes; le deuxième c'est d'offrir une règle encore plus simple pour la formation de ces formules. En effet, il est facile de remarquer[21]

1° que le premier terme de l'une quelconque de ces équations, est formé du premier membre de l'équation précédente, multiplié par la première des lettres qu'elle ne renferme point, cette lettre étant affectée de la marque qui suit immédiatement la plus haute de celles qui entrent dans ce même membre ;

---

[21] Pour éclaircir ce qui va suivre dans les 1°, 2°, etc., voici littéralement cette formation pour *n* = 1, 2, et 3 :
- Pour une équation à une inconnue, *ax* = 0, la condition pour des solutions non nulles est *a* = 0 ;
- Pour 2 équations à 2 inconnues, *ax+by* = 0, *a'x+b'y* = 0, le 1° donne *ab'* et le 2° donne –*a'b*, donc la condition pour des solutions non nulles est bien celle que Bézout a annoncée plus haut     *ab'-a'b* = 0 ;
- Pour 3 équations à 3 inconnues (notations analogues), le 1° donne *(ab'-a'b)c″*, le 2° donne –*(ab″-a″b)c'*, et le 3° donne –*(a″b-a'b″)c* , donc la condition pour des solutions non nulles est, là aussi, celle annoncée au dessus *(ab'-a'b)c″+(a″b-ab″)c'+(a'b″-a″b')c* = 0.

2° Le deuxième terme se forme du premier, en changeant dans celui-ci la plus haute marque en celle qui est immédiatement au–dessous & réciproquement, & de plus en changeant les signes ;

3° Le troisième, se forme du premier, en changeant dans celui-ci la plus haute marque en celle de deux numéros au–dessous & réciproquement, & de plus en changeant les signes ;

4° Le quatrième, se forme du premier, en changeant dans celui-ci la plus haute marque en celle de trois numéros au–dessous & réciproquement, & changeant les signes, & toujours de même pour les suivans. » [*Ibid*. p.294]

Il ajoute en corollaire : « Chacun des termes de l'équation de condition a donc essentiellement le même nombre de facteurs, & ces facteurs sont tellement combinés que jamais, dans un même terme, il ne s'y en rencontre deux qui appartiennent à une même inconnue » [*Ibid*.].

Sa règle d'écriture a deux caractéristiques, elle donne en même temps les termes et leurs signes et elle procède par induction sur le nombre d'inconnues et donc de coefficients d'une équation. Il est remarquable, en plus de sa simplicité et des caractéristiques précédentes, qu'elle donne les mises en forme suivantes :

1°) $ab'c''-ac'b''+ca'b''-ba'c''+bc'a''-cb'a'' = (ab'-a'b)c''+(a''b-ab'')c'+(a'b''-a''b')c$

2°)  $ab'c''d'''-ab'd''c'''+ad'b''c'''-da'b''c'''-ac'b''d'''+ac'd''b'''-ad'c''b'''+da'c''b'''+ca'b''d'''-$
$ca'd''b'''+cd'a''b'''-dc'a''b'''-ba'c''d'''+ba'd''c'''-bd'a''c'''+db'a''c'''+bc'a''d'''-bc'd''a'''+bd'c''a'''-$
$db'c''a'''-cb'a''d'''+cb'd''a'''-cd'b''a'''+dc'b''a'''$  $= [(ab'-a'b)c''+(a''b-ab'')c'+(a'b''-a''b')c]d'''$ + $[(a'b-ab')c'''+(ab'''-a'''b)c'+(a'''b'-a'b''')c]d''$ + $[(a'''b-ab''')c''+(ab''-a''b)c'''+(a''b'-a'''b'')c]d'$ + $[(a'b'''-a'''b')c''+(a''''b''-a''b''')c'+(a'''b'-a'b''')c'']$ $d$. [voir *Ibid*. p. 293]

Bien que le contexte théorique et algébrique soit très différent, on reconnaît ici ce qui fut appelé postérieurement, le développement d'un déterminant suivant les éléments d'une ligne ou d'une colonne, faisant apparaître les déterminants mineurs correspondants.

Bézout revient ensuite au calcul de la résultante. Il constate que les méthodes d'Euler et de Cramer, très performantes de ce point de vue, pour deux équations à deux inconnues, ne le sont plus dès qu'on dépasse ce nombre d'équations car « telle est la nature de ces méthodes, qu'elle exige que pour éliminer on compare les équations deux à deux : or [...] ce procédé conduit à des équations beaucoup plus élevées qu'il ne faut. » [*Ibid.*]. En effet, les méthodes de ses prédécesseurs[22], par exemple dans le cas de trois équations $(A)$, $(B)$, $(C)$ à trois inconnues $x, y, z$, consistaient à éliminer d'abord $x$ entre $(A)$ et $(B)$, ce qui donnait une nouvelle équation $(D)$ d'inconnues $y$ et $z$, puis $x$ entre $(B)$ et $(C)$, pour obtenir une équation $(E)$ en $y, z$, aussi. Ils éliminaient alors $y$ entre $(D)$ et $(E)$ et arrivaient à une équation en $z$, la résultante de (A), (B), (C), qui, à cause de tous les calculs intermédiaires, contenait presque toujours des facteurs superflus et des racines étrangères aux équations de départ.

Pour remédier à ces défauts, Bézout entreprend, pour éliminer une inconnue :

- de traiter toutes les équations en même temps ;

- de ramener l'élimination à la résolution d'un système linéaire ;.

- enfin de trouver *a priori* le degré de la résultante – c'est-à-dire du polynôme en la ou les inconnue(s) restante(s).

Étant données $n$ équations à $n$ inconnues de degrés quelconques, sa démarche consiste

- à multiplier chacune des équations par un polynôme à coefficients indéterminés ;

---

[22] Les méthodes de ses prédécesseurs sont celles que donne d'Alembert :
« Quand il y a plus de deux inconnues, par exemple, *x, y, z, &c.* on réduit d'abord les inconnues à une de moins ; on fait évanoüir *x* ou *y*, &c. en traitant *z* & les autres comme une constante ; ensuite on réduit les inconnues restantes à une de moins, & ainsi du reste. Cela n'a aucune difficulté. » [*Encyclopédie* t. V, 1756, article « Évanouir »]. Sa conclusion péremptoire (« Cela n'a aucune difficulté ») a peut-être fait sourire Bézout quand il en a pris connaissance.

- à égaler la somme de tous ces produits à zéro, obtenant ainsi ce qu'il appelle « l'équation-somme » ;

- à annuler dans cette équation-somme tous les coefficients de l'inconnue qu'il a choisie d'éliminer, grâce aux indéterminées des différents polynômes multiplicateurs.

Il aborde ainsi simultanément deux problèmes : déterminer les degrés des polynômes multiplicateurs et celui de la résultante (partant du principe que, comme pour le cas de deux équations à deux inconnues, il n'y a, les inconnues restantes étant choisies, qu'une résultante à un coefficient multiplicateur près).

Étienne Bézout calcule d'abord le degré de la résultante pour deux équations à deux inconnues :

« non que je prétende, justifie-t-il, décider ma méthode préférable à celle que M$^{rs}$. Euler & Cramer ont donnée pour ce cas seulement, mais parce que cette méthode étant uniforme, j'ai cru me rendre plus clair en fortifiant l'analogie par la réunion de ce cas avec les autres, & en même temps parce que dans un travail aussi long que l'est souvent celui de l'élimination, il n'est pas inutile de multiplier les méthodes sur lesquelles les Calculateurs peuvent porter leur choix » [*Ibid*.p. 291].

$$\text{Soit les deux équations} \begin{cases} Ax^m + Bx^{m-1} + Cx^{m-2} + Dx^{m-3} + Ex^{m-4} + \dots\dots\dots + V = 0 \\ A'x^{m'} + B'x^{m'-1} + C'x^{m'-2} + D'x^{m'-3} + E'x^{m'-4} + \dots\dots + V' = 0 \end{cases},$$

dans lesquelles $A$, $B$, $C$, $D$, $E$, etc., $A'$, $B'$, $C'$, $D'$, $E'$, etc., sont des polynômes en la deuxième inconnue $y$, de degrés respectifs : $p$, $p+1$, $p+2$, $p+3$, etc., et $p'$, $p'+1$, $p'+2$, $p'+3$, etc.

Il multiplie les deux équations respectivement, par les polynômes $(L)$ et $(L')$ suivants :

$(L) \quad Mx^n + Nx^{n-1} + Px^{n-2} + Qx^{n-3} + Rx^{n-4} + \dots\dots\dots\dots + T$

$(L') \quad M'x^{n'} + N'x^{n'-1} + P'x^{n'-2} + Q'x^{n'-3} + R'x^{n'-4} + \dots\dots\dots + T'$

et obtient donc, pour que chaque puissance de x disparaisse dans l'équation-somme, le système ci-dessous :

$$
\left\{
\begin{array}{l}
\qquad\qquad\qquad\qquad\qquad\qquad AM \;+\; A'M' \;=\; 0\\
\qquad\qquad\qquad\qquad AN \;+\; A'N' \;+\; BM \;+\; B'M' \;=\; 0\\
\qquad\qquad AP \;+\; A'P' \;+\; BN \;+\; B'N' \;+\; CM \;+\; C'M' \;=\; 0\\
AQ \;+\; A'Q' \;+\; BP \;+\; B'P' \;+\; CN \;+\; C'N' \;+\; DM \;+\; D'M' \;=\; 0\\
\qquad\qquad\quad \& c.\\
\qquad\qquad VT \;+\; V'T' \;=\; 0
\end{array}
\right.
$$

avec les conditions $m+n=m'+n'$, et $m+n+1=n+1+n'+1$ (*i.e.* autant d'équations que d'inconnues). Les valeurs de $n$ et $n'$ sont alors fixées, $n = m'-1$ et $n' = m-1$, et le déterminant du système doit être nul car sinon les inconnues $M$, $M'$, $N$, $N'$, etc., devraient toutes l'être et cela impliquerait que les deux équations de départ n'ont pas de racines communes.

Voilà donc la nouveauté, précédemment annoncée, de l'approche de Bézout : il définit la résultante du système de départ - de deux équations à deux inconnues - comme l'équation obtenue en écrivant que le déterminant du système linéaire induit est égal à zéro.

Il remarque que les coefficients d'une même inconnue sont des polynômes dont les degrés sont en progression arithmétique de raison 1. Il fait alors appel au lemme II qu'il a démontré au début de son mémoire :

« Si on a un nombre quelconque $n$ de quantités *a, b, c, d, e*, &c. & qu'au dessous de chacune de ces quantités, on écrive les progressions arithmétiques suivantes,

| | | | | | |
|---|---|---|---|---|---|
| *a* | *b* | *c* | *d* | *e* | &c. |
| *a+k* | *b+k* | *c+k* | *d+k* | *e+k* | &c. |
| *a+2k* | *b+2k* | *c+2k* | *d+2k* | *e+2k* | &c. |
| *a+3k* | *b+3k* | *c+3k* | *d+3k* | *e+3k* | &c. |
| &c. | &c. | &c. | &c. | &c. | &c. |

en continuant les progressions jusqu'à ce que le nombre des termes de chacune soit égal au nombre $n$ des quantités, je dis que dans quelque ordre qu'on ajoute $n$ termes de cette progression, pourvu qu'on n'y comprenne jamais deux termes d'une même colonne ni deux termes d'une même bande, la somme sera toujours la même $\& = S + \dfrac{kn(n-1)}{2}$, $S$ marquant la somme des termes qui composent la première bande » [Bézout 1767c, p. 296]

Il en est de même s'il y a des lacunes dans le tableau, en sommant comme premiers termes, ceux qui sont ou devraient être en première ligne.

Par application immédiate de ce lemme aux suites arithmétiques des degrés en $y$ des coefficients $A$, $A'$, $B$, $B'$, $C$, $C'$, etc., qui se retrouvent dans le système linéaire en $M$, $M'$, $N$, $N'$, $P$, $P'$, etc., c'est-à-dire

| $p$ | $p'$ | | | | |
|---|---|---|---|---|---|
| $p+1$ | $p'+1$ | $p$ | $p'$ | | |
| $p+2$ | $p'+2$ | $p+1$ | $p'+1$ | $p$ | $p'$ |
| $p+3$ | $p'+3$ | $p+2$ | $p'+2$ | $p+1$ | $p'+1$ |
| etc. | etc. | etc. | etc. | etc. | etc. |

on trouve que le degré de chaque terme du déterminant du système sera au plus égal à

$$(2p-n)\frac{n+1}{2} + (2p'-n')\frac{n'+1}{2} + \frac{m+n+1}{2}(m+n).$$

Or, l'annulation de ce déterminant donne l'équation résultante du système de départ donc en appelant G le maximum possible du degré de cette résultante, on a,

$$G = (2p-n)\frac{n+1}{2} + (2p'-n')\frac{n'+1}{2} + \frac{m+n+1}{2}(m+n), \text{ et comme } n = m'\text{-}1 \text{ et } n' = m\text{-}1$$

$$\mathbf{G = mm' + mp' + m'p}.$$

On peut remarquer plusieurs autres résultats intéressants qui ressortent de ce texte :

- deux polynômes $P(x,y)$ et $Q(x,y)$ étant donnés d'ordre respectivement $m$ et $m$' en $x$, Bézout a, en fait, montré qu'il existe un polynôme non nul $L_1(x,y)$ d'ordre au plus $m$'-1 en $x$, et un polynôme non nul $L_2(x,y)$ d'ordre au plus $m$-1 en $x$, tels que $L_1.P + L_2.Q = 0$ soit une équation en $y$, si et seulement si, le déterminant du système obtenu en égalant à zéro tous les coefficients des différentes puissances de $x$, dans $L_1.P + L_2.Q$, est nul ;

- la résultante qu'il obtient, étant le déterminant d'un système dont les coefficients sont des polynômes en $y$, Bézout est sûr que c'est aussi un polynôme en $y$, puisqu'on ne fait que des multiplications et des sommes de polynômes. Euler et Cramer manipulant des fractions et des racines de polynômes, n'ont pas montré ce point de façon convaincante.

- si les deux équations de départ sont de degrés respectifs $m$ et $m'$, alors $p = p' = 0$, et on retrouve, dans sa formule, pour maximum du degré de la résultante le produit $m.m'$ des degrés des équations.

Par contre, Bézout a oublié de préciser que les deux courbes définies par les équations du système de départ ne doivent pas avoir de composante commune, car sinon il y aurait une infinité de points communs, ceux de la composante.

Donc, à cela près, Bézout a démontré rigoureusement que deux courbes planes de degrés respectifs $m$ et $m'$, se coupent au plus en $m.m'$ points. Avant même la généralisation qu'il fera en 1779, pour un nombre quelconque d'équations et d'inconnues, cela justifie déjà que le théorème correspondant porte son nom.

Par contre, Euler a montré que sa résultante ne contient aucun facteur superflu et donc ne peut donner aucune racine étrangère au problème, alors que, si Bézout donne une condition nécessaire pour que le système ait des racines, il ne démontre pas qu'elle est aussi suffisante[23]. Étant parvenu à une équation de même degré qu'Euler et Cramer, Bézout a sans doute estimé qu'il ne pouvait y avoir de facteur superflu dans sa résultante, d'après le résultat même de ses prédécesseurs.

Après cet exposé de son travail, il convient de revenir sur les caractéristiques de la méthode de Bézout. Celui-ci utilise, au départ, le principe des coefficients indéterminés tel qu'il est exposé par Euler dans son *Introductio in analysin infinitorum*, [Euler 1748, t. 2, ch. XIX], pour trouver la résultante de deux équations à deux inconnues. Mais Euler n'utilise pas la méthode des coefficients indéterminés pour obtenir, *a priori*, le degré de la résultante, alors que c'est ainsi que Bézout l'utilise. Celui-ci a donc raison de dire qu'il ne suit pas la même méthode qu'Euler et Cramer, *pour calculer le degré de la résultante* [souligné par nous].

---

[23] La démonstration de la condition suffisante est aujourd'hui aisée, grâce au Résultant.

Par ailleurs, il faut souligner la différence entre les méthodes d'Euler [1748] et de Bézout [1767c] pour l'obtention de l'équation résultante. Certes, elles suivent la même idée de départ : multiplier chacune des équations par un polynôme à coefficients indéterminés, chacun de ces deux polynômes ayant pour degré, par rapport à l'inconnue qui apparaît dans l'écriture, celui de cette inconnue dans l'équation qu'il ne multiplie pas, moins un.

À partir de là, Euler identifie les coefficients des deux polynômes produits, ce qui lui permet d'obtenir, $m$ étant le degré en $y$ de la première équation et $n$ celui de la seconde, un système d'équations linéaires de $m+n-1$ équations à $m+n-2$ inconnues. Mais chez lui, à cause de certains coefficients choisis[24], le système n'est pas homogène. Il calcule les coefficients indéterminés avec les $m+n-2$ premières équations, et, en plaçant ces valeurs dans la dernière équation, il obtient la résultante.

Quels sont les défauts de ce procédé ? D'une part, le système n'étant pas homogène, il peut ne pas avoir de solutions. D'autre part, lorsque les solutions existent, ce sont souvent des fractions rationnelles et non des polynômes. Enfin, Euler utilise la méthode des coefficients indéterminés de façon classique, c'est-à-dire en identifiant des termes et en calculant les coefficients pour les remplacer dans une des équations obtenues. Il suppose implicitement que ces coefficients existent, et *il cherche leurs valeurs, valeurs dont il a besoin* [souligné par nous].

Étienne Bézout, après avoir multiplié lui aussi chaque équation (d'ordres respectifs $m$ et $m'$ en $x$) par un polynôme à coefficients indéterminés (d'ordres respectifs $m'-1$ et $m-1$ en $x$), fait la somme des deux produits et égale à zéro les coefficients de tous les termes en *x*. Mais chez lui, on l'a vu, le système obtenu de *m+m'* équations à *m+m'* inconnues (les coefficients indéterminés), est homogène. Il a donc toujours des solutions et l'équation résultante est la condition d'existence de solutions non toutes nulles, c'est-à-dire, l'annulation du déterminant. Ce dernier étant constitué de sommes de produits de polynômes, est lui-même un polynôme. Bézout ne travaille ainsi que sur des polynômes. On le voit, la même idée de départ est utilisée de façon plus élégante et surtout plus sûre chez Bézout.

---

[24] Dans chacun des polynômes multiplicateurs à coefficients indéterminés, un des coefficients est déterminé. Le coefficient du terme en $y^{n-1}$ du premier polynôme multiplicateur est celui du terme en $y^n$ de la seconde équation et réciproquement, le coefficient du terme en $y^{m-1}$ du second polynôme multiplicateur est celui du terme en $y^m$ de la première équation. Il y a donc des termes connus, non nuls, dans les coefficients de $y$ de chaque produit.

D'autre part, Bézout, lui, se sert seulement des coefficients indéterminés comme outils pour construire sa résultante. Il n'a pas besoin de connaître leurs valeurs et d'ailleurs il ne les cherche pas. Il veut uniquement connaître leur condition d'existence, qui est l'équation cherchée. Il faut là remarquer l'originalité de l'application que fait Bézout des coefficients indéterminés. Ce n'est pas du tout l'utilisation habituelle, et nous verrons par la suite que, dans son ouvrage de 1779, il continuera à les employer d'une façon constructive qui lui est propre [Alfonsi 2006].

Pour illustrer ce qui précède, voici les deux méthodes appliquées à un exemple d'Euler,

$$\begin{cases} Py^2 + Qy + R = 0 \\ py^3 + qy^2 + ry + s = 0 \end{cases}$$ , où *P, Q, R, p, q, r, s*, sont des polynômes en *x* :

-<u>Méthode (et notations) d'Euler</u> [Euler 1748, t. 2, ch. XIX] :

Il multiplie la première équation par $py^2 + ay + b$ et la seconde par $Py + A$, où seuls *a, b, A*, sont les coefficients indéterminés. Il obtient, en identifiant les coefficients d'une même puissance de *y* pour chacun des deux produits obtenus, le système linéaire non homogène de quatre équations à trois inconnues *a*, *A*, et *b* :

$$Pa + Qp = pA + qP = \alpha \qquad\qquad Pb + Qa + Rp = qA + rP = \beta$$

$$Qb + Ra = rA + sP \qquad\qquad Rb = sA$$

Euler tire de la première équation, $a = \dfrac{\alpha - Qp}{P}$ et $A = \dfrac{\alpha - qP}{p}$, puis de la deuxième,

$$b = \frac{\beta - Qa - Rp}{P} = \frac{\beta - Rp}{P} - \frac{\alpha Q - Q^2 p}{P^2} \qquad \text{et} \qquad \beta = rP + \frac{\alpha q - q^2 P}{p}$$

Substituant cette valeur de β dans $b$, il obtient $b$ en fonction de α uniquement, comme c'est déjà le cas pour $a$ et $A$. En substituant maintenant, ces valeurs de $b$, de $a$ et de $A$ dans la troisième équation, il obtient α, puis enfin, $a$, $A$ et $b$, en fonction de $P$, $Q$, $R$, $p$, $q$, $r$, $s$. La quatrième équation lui donne alors la résultante en $x$.

- <u>Méthode de Bézout sur cet exemple</u>[25] :

Il multiplie la première équation par $My^2 + Ny + L$, et la deuxième par $M'y + N'$. Il écrit l'« équation-somme », $(My^2 + Ny + L)(Py^2 + Qy + R) + (M'y + N')(py^3 + qy^2 + ry + s) = 0$.

Il égale à 0 chaque coefficient des puissances de $y$, ce qui lui donne un système linéaire homogène de cinq équations à cinq inconnues $M$, $M'$, $N$, $N'$, $L$ :

$$\begin{cases} & & & & MP & + & M'p & = & 0 \\ & & NP & + & N'p & + & MQ & + & M'q & = & 0 \\ LP & + & NQ & + & N'q & + & MR & + & M'r & = & 0 \\ LQ & + & NR & + & N'r & & & + & M's & = & 0 \\ LR & & & + & N's & & & & & = & 0 \end{cases}$$

L'équation qui égale à zéro le déterminant de ce système, donne la résultante en $x$ cherchée.

En revanche, si l'on examine les deux méthodes ([Euler 1748] et [Bézout 1767c]) du point de vue de l'algèbre linéaire actuel, on constate qu'elles nous conduisent à la même matrice carrée que nous noterons **R** (matrice résultante), qui est d'ordre $m+m'$, et qui avec les

notations polynomiales de Bézout est égale à **R** =
$$\begin{pmatrix} A & 0 & 0 & ... & A' & 0 & 0 & ... \\ B & A & 0 & ... & B' & A' & 0 & ... \\ C & B & A & ... & C' & B' & A' & ... \\ D & C & B & ... & D' & C' & B' & ... \\ etc. & D & C & ... & etc. & D' & C' & ... \\ 0 & etc. & D & ... & 0 & etc. & D' & ... \\ 0 & 0 & etc. & ... & 0 & 0 & etc. & ... \\ 0 & 0 & 0 & ... & 0 & 0 & 0 & ... \end{pmatrix}$$

---

[25] Sachant que Bézout ne l'a pas traité, ni d'ailleurs aucun autre dans son mémoire de 1764, pour illustrer sa méthode générale de calcul du degré de la résultante

C'est le déterminant de **R** qu'on appelle couramment aujourd'hui le *Résultant*[26] du système

$$\begin{cases} Ax^m + Bx^{m-1} + Cx^{m-2} + Dx^{m-3} + Ex^{m-4} + \ldots\ldots\ldots\ldots + V = 0 \\ A'x^{m'} + B'x^{m'-1} + C'x^{m'-2} + D'x^{m'-3} + E'x^{m'-4} + \ldots\ldots + V' = 0 \end{cases}$$

et dont l'égalité ou non à zéro permet de savoir si le système ci-dessus a des solutions.

Si on peut attribuer à Euler l'idée qui a conduit, bien plus tard, à cette matrice résultante, c'est Bézout qui a le premier pensé à utiliser l'annulation de son déterminant comme condition nécessaire à la résolution du système des deux équations. Pour deux équations à deux inconnues la détermination du degré est donc pleinement réussie, car les conditions permettent de déterminer sans ambiguïté *n*, *n' et G*. Il n'en sera pas de même quand Bézout augmentera les nombres d'équations et d'inconnues considérées. Pour trois équations à trois inconnues, *x, y, z,* il est obligé, puisqu'il n'envisage que l'élimination d'une inconnue dans son calcul, d'appliquer sa méthode deux fois, c'est à dire qu'il devra obtenir trois polynômes pour éliminer *x* et obtenir un premier polynôme en *y, z* ; puis trois autres pour un second polynôme en *y, z* ; enfin, il devra éliminer *y* ou *z* dans les deux derniers polynômes trouvés.

Il commence donc par éliminer *x*. Les notations sont analogues aux précédentes : Soit trois équations de degrés respectifs en *x : m, m', m''*, dans lesquelles les coefficients *A, B, C,* etc., *A', B', C', etc., A'', B''*, etc., sont des polynômes à deux inconnues y et z, telles que les degrés de *A, B, C*, etc., *A', B', C',* etc., *A'', B'', C''*, etc., soient respectivement *p, p+1, p+2,* etc., *p', p'+1, p+2*, etc., *p'', p''+1, p''+2*, etc. Les polynômes multiplicateurs de chaque équation ont pour degrés respectifs en *x, n, n', n''*. Bézout suppose, pour que les termes de plus haut degré de l'équation somme puissent s'annuler, que $m+n = m'+n'$ et $m+n \geq m''+n''$. On doit aussi avoir $m+n+1 = n+1+n'+1+n''+1$, (autant d'équations que d'inconnues), donc $n'=m-n''-2$ *et* $n=m'-n''-2$.

---

Par un raisonnement analogue au précédent, sur le degré du déterminant du système obtenu en annulant les coefficients de $x$ dans l'équation-somme d'ordre $n$, il obtient, $G$ étant le degré de ce premier polynôme en $y, z$ :

$$G = mm' + pm' + p'm - m - m' + m'' - p - p' + p'' + 1 - (p + p' - p'' + m + m' - m'' - n'' - 2)n'',$$

où $n''$ est indéterminé, et pour avoir G minimum, il annule la dérivée de G par rapport à $n''$, ce qui donne $n'' = \dfrac{(m + m' + m'' + p + p' + p'')}{2} - m'' - p'' - 1$ .

Cette valeur de $n''$ ne pourra pas toujours être prise car on doit avoir à la fois : $n''$ positif, le terme entre parenthèses pair, $m + n \geq m'' + n''$ et $m' - 2 \geq n''$, puisqu'on suppose $m' \leq m$. Bézout choisit alors $n'' = [(m + m' + m'' + p + p' + p'' - \boldsymbol{\alpha})/2] - m'' - p'' - 1$ où $\alpha$ est « la plus petite valeur possible ». Il ne trouve pas, dans ce cas, de formule qui lui donne de façon certaine le meilleur résultat pour G, puisque le choix de $n''$ reste aléatoire. Il reconnaît le caractère insatisfaisant de son résultat puisqu'il écrit :

« Si l'on ne peut parvenir à donner à $n, n', n''$, des valeurs positives, qu'en rendant G plus grand qu'il ne serait par la combinaison des équations deux à deux, on aura recours à ce dernier moyen » [Bézout 1767c, p. 305], ce qui le ramène alors à la méthode de ses prédécesseurs.

## 3.     La méthode du Bézoutien[27]

Il revient maintenant sur le cas de deux équations, mais seulement pour calculer la résultante d'une autre façon, puisque, dorénavant il connaît le degré qu'elle doit avoir. Il envisage d'abord le cas de deux équations à deux inconnues de même degré $m$ en $x$.

---

[27] La matrice **B** que nous considèrerons dans ce paragraphe, a été appelée par J.J. Sylvester [1853], le « **Bezoutiant** », reconnaissant ainsi le rôle de Bézout dans le processus qui l'a faite naître, et c'est le nom que cette matrice a gardé (on trouve aussi en anglais, « Bezoutian », ou « Bezout matrix »). En français l'appellation de Sylvester a été transformée en Bézoutien. C'est ce dernier nom que l'on trouve le plus souvent dans la littérature française contemporaine [Levavasseur 1907], [Elkadi 2007], et c'est donc celui que nous avons choisi d'employer dans cet article. Sur le Bézoutien et ses propriétés, en plus des ouvrages déjà cités, voir [Wimmer 1990], [Fuhrmann et Helmke 1989], [Fuhrmann 1996], [Lerer et Haimovici 1995], [Lerer et Rodman 1996, 1999], par exemple.

Contrairement à la précédente, cette méthode d'obtention de la résultante est entièrement une création de Bézout. On verra qu'elle n'utilise pas de coefficients indéterminés, mais seulement les coefficients des équations de départ. De plus le calcul sera plus facile que celui du résultant puisque l'ordre du déterminant obtenu sera, nous le verrons, seulement $m$, au lieu de $m+m' = 2m$.

Il considère le système

$$\begin{cases} Ax^m + Bx^{m-1} + Cx^{m-2} + Dx^{m-3} + Ex^{m-4} + \ldots\ldots\ldots + Ux + V = 0 \\ A'x^m + B'x^{m-1} + C'x^{m-2} + D'x^{m-3} + E'x^{m-4} + \ldots\ldots + U'x + V' = 0 \end{cases}$$

et il multiplie :

- la première équation par $A'$ et la deuxième par $A$, la différence est de degré $m$-1 ;

- la première équation par $A'x+B'$ et la deuxième par $Ax+B$, la différence est de degré $m$-1 ;

- la première équation par $A'x^2+B'x+C'$ et la deuxième par $Ax^2+Bx+C$, la différence est de degré $m$-1 ; etc.

jusqu'à avoir $m$ équations, chacune de degré $m$-1.

Il considère chaque puissance de $x$ comme inconnue, et obtient donc un système de $m$ équations à $m$-1 inconnues. Grâce au résultat sur les équations linéaires donné au début, il peut dire que le déterminant du système de $m$ équations à $m$ inconnues, obtenu à partir de celui qu'il a construit, en multipliant la colonne des constantes par une même lettre, doit être égal à zéro pour qu'il y ait des solutions[28].

Notons les deux équations de départ $f(x)=0$ et $g(x)=0$, et les polynômes multiplicateurs, à partir de $f$ : $f_m = A$ , $f_{m-1} = Ax+B$ , $f_{m-2} = Ax^2+Bx+C$ ,…etc…, $f_{m-i} = Ax^i+Bx^{i-1}+$ etc., etc.,

---

[28] L'exemple donné *supra* au début du III.1, $\begin{cases} P(x, y) = A(y)x^2 + B(y)x + C(y) = 0 \\ Q(x, y) = A'(y)x^2 + B'(y)x + C'(y) = 0 \end{cases}$ est traité par Bézout de cette façon. Il calcule A'P – AP', puis (A'x + B')P – (Ax + B)P', et obtient ainsi $\begin{cases} [A'(y)B(y) - A(y)B'(y)]x + A'(y)C(y) - A(y)C'(y) = 0 \\ [A'(y)C(y) - A(y)C'(y)]x + B'(y)C(y) - B(y)C'(y) = 0 \end{cases}$ d'où la résultante donnée précédemment $[A(y)B'(y) - B(y)A'(y)] [C'(y)B(y) - B'(y)C(y)] = [A'(y)C(y) - A(y)C'(y)]^2$

$f_1 = Ax^{m-1} + Bx^{m-2} + Cx^{m-3} + \dots + U$,   et à partir de $g$ :   $g_m = A'$,   $g_{m-1} = A'x + B'$,   $g_{m-2} = A'x^2 + B'x + C'$,

etc.,   $g_{m-i} = A'x^i + B'x^{i-1} + \text{etc.}$,   etc.,   $g_1 = A'x^{m-1} + B'x^{m-2} + C'x^{m-3} + \dots + U'$,

on peut écrire, avec nos notations actuelles, le système obtenu par Bézout sous la forme

$$B . \begin{pmatrix} x^{m-1} \\ x^{m-2} \\ etc. \\ \dots \\ x \\ 1 \end{pmatrix} = f(x) \begin{pmatrix} g_m \\ g_{m-1} \\ etc. \\ \dots \\ g_2 \\ g_1 \end{pmatrix} - g(x) \begin{pmatrix} f_m \\ f_{m-1} \\ etc. \\ \dots \\ f_2 \\ f_1 \end{pmatrix} \qquad \text{où la matrice } \mathbf{B} \text{ est égale à}$$

$$B = \begin{pmatrix} A & 0 & 0 & \dots & 0 \\ B & A & 0 & \dots & 0 \\ C & B & A & \dots & 0 \\ \dots & \dots & \dots & \dots & \dots \\ U & T & S & \dots & A \end{pmatrix} \begin{pmatrix} B' & C' & D' & \dots & V' \\ C' & D' & \dots & \dots & 0 \\ D' & \dots & V' & \dots & 0 \\ \dots & \dots & \dots & \dots & \dots \\ V' & 0 & 0 & \dots & 0 \end{pmatrix} - \begin{pmatrix} A' & 0 & 0 & \dots & 0 \\ B' & A' & 0 & \dots & 0 \\ C' & B' & A' & \dots & 0 \\ \dots & \dots & \dots & \dots & \dots \\ U' & T' & S' & \dots & A \end{pmatrix} \begin{pmatrix} B & C & D & \dots & V \\ C & D & \dots & \dots & 0 \\ D & \dots & V & \dots & 0 \\ \dots & \dots & \dots & \dots & \dots \\ V & 0 & 0 & \dots & 0 \end{pmatrix}$$

L'équation obtenue en écrivant l'annulation du déterminant de **B**, (le Bézoutien[29]), est, d'après ce qui précède, la résultante du système de départ, donc une condition nécessaire et suffisante pour que ce système ait des solutions.

En fait, pour une variable, le système a une solution si et seulement si Dét **B** = 0 ou si et seulement si Dét **R** = 0 [Sylvester 1840, 1853]

Il envisage maintenant le cas où les deux équations à 2 inconnues sont de degrés différents en $x$, $m$ et $m'$, avec $m' < m$ :

$$\begin{cases} Ax^m + Bx^{m-1} + Cx^{m-2} + Dx^{m-3} + Ex^{m-4} + \dots\dots\dots + V = 0 \\ A'x^{m'} + B'x^{m'-1} + C'x^{m'-2} + D'x^{m'-3} + E'x^{m'-4} + \dots\dots + V' = 0 \end{cases}$$

Bézout essaie le même genre de procédé, mais cela s'avère moins simple. Il multiplie

- la première par $A'$ et la deuxième par $A\,x^{m-m'}$ et fait la différence qui est de degré $m$-1 en $x$ ;

- la première par $A'x + B'$ et la deuxième par $Ax^{m-m'+1} + Bx^{m-m'}$, la différence est de degré $m$-1 en $x$ ;

---

- la première par $A'x^2+B'x+C'$ et la deuxième par $Ax^{m-m'+2} + Bx^{m-m'+1} + Cx^{m-m'}$, la différence est de degré $m$-1 en $x$ ; etc., jusqu'à avoir $m'$ équations, chacune de degré $m$-1.

Puis il multiplie chacune de ces équations par un coefficient indéterminé, et la deuxième équation de départ par $M'x^{m-m'-1} + N'x^{m-m'-2} + P'x^{m-m'-3} + \ldots\ldots + T$ ; une fois tous ces produits obtenus, il les ajoute et égale à 0 la somme des coefficients de chaque puissance de $x$. Par les formules du lemme II, il obtient l'équation de condition.

Dans ce cas, il retombe donc sur la première méthode[30] puisque ce n'est plus un système d'équations où les inconnues sont des puissances de $x$ qu'il obtient, mais un système où les inconnues sont les coefficients indéterminés introduits. Bézout poursuit par des exemples et recherche des procédés pour simplifier les calculs, mais ne trouve plus un niveau de généralité et de simplicité intéressant.

Il le trouve encore moins lorsqu'il veut passer à des équations à un plus grand nombre d'inconnues. Il est, comme dans sa première méthode, confronté aux mêmes difficultés de détermination des exposants et à des calculs longs et aléatoires. Il en est d'ailleurs très conscient : « Au reste, je crois ces méthodes encore très susceptibles de perfection, & il y a un grand nombre de cas où en suivant les principes sur lesquels elles sont fondées, on parvient à trouver des routes plus faciles. [...] c'est un travail auquel j'invite ceux qui seront assez heureux pour avoir plus de temps à dépenser que moi. » [Bézout 1767c, p. 329]

Cette dernière remarque figurant dans la dernière partie du mémoire lu le 29 février 1764, fait sans doute allusion à la perspective de nouvelles charges liées à la formation des gardes du Pavillon et de la Marine, le conduisant à se rendre à Brest dès mars 1764.

Ce point est important car il explique en partie le fait que Bézout se soit exclusivement consacré (le manque de temps nécessaire à la recherche l'obligeant à faire des choix) à l'analyse algébrique finie (plus spécialement l'élimination) pendant le reste de sa vie et il explique aussi son peu de productions en recherche mathématique jusqu'en 1779. Nous y reviendrons plus loin.

---

[30] On verra que c'est en rédigeant, en 1766, la partie algèbre de son cours pour la Marine que Bézout trouvera l'application du Bézoutien au cas où m et m' sont différents. C'est d'ailleurs dans son cours qu'il exposera cette méthode, nous y reviendrons.

## 4. L'« Identité de Bézout »

Dans le tout dernier paragraphe Bézout écrit :

« On peut encore appliquer à beaucoup d'autres usages cette manière de combiner les équations, particulièrement à la recherche du commun diviseur de plusieurs quantités complexes. En effet, si deux, trois ou un plus grand nombre de quantités complexes ont un diviseur commun composé, par exemple, de $x$, & de telles autres quantités qu'on voudra, on peut supposer que chacune de ces quantités est zéro, parce qu'elle le deviendroit en effet, si on mettoit pour $x$ sa valeur qu'on auroit en égalant ce diviseur à zéro : alors ou le diviseur est d'une, ou de deux ou d'un plus grand nombre de dimensions ; il n'y a donc qu'à chercher quels sont les polinomes indéterminés par lesquels il faudroit multiplier chaque équation pour qu'en égalant à zéro la somme des produits, il n'y restât que les deux derniers termes, si le commun diviseur doit être d'une dimension, ou les trois derniers, s'il doit être de deux dimensions, & ainsi de suite. Je me contente d'indiquer cet usage. » [Bézout 1767c, p. 337-338]

L'énoncé peut laisser penser qu'il envisage plusieurs inconnues alors qu'il n'en considère qu'une, celle qu'il cite $x$, et il énonce donc la propriété pour une seule variable. On peut comprendre ce qui suit : Bézout considère la première méthode qu'il a exposée, celle de la multiplication par des polynômes indéterminés ; il se place dans le cas où le résultant est égal à zéro, nécessaire pour que les polynômes aient des solutions communes ; dans ce cas, il essaye d'annuler seulement les coefficients des termes en $x$ pour des exposants de $x$ supérieurs à 1, ou 2, etc., ce qui va lui donner $m+n$-2, ou $m+n$-3, etc.. équations à $m+n$ inconnues (les coefficients indéterminés), qu'il pourra toujours résoudre en fonction de deux, ou trois, etc., coefficients choisis ; grâce à cela, il aura avec les deux (ou trois, etc.) derniers coefficients des puissances de $x$, deux (ou trois, etc.) équations à deux (ou trois, etc.) inconnues ; suivant la compatibilité de ce dernier système, il pourra conclure s'il a ou non atteint le degré du PGCD.

On peut constater que Bézout considère le PGCD de plusieurs polynômes comme la somme des produits de chacun de ces polynômes par un polynôme à coefficients indéterminés. On retrouvera cette idée en 1853 chez Sylvester[31] puis avec l'émergence de la notion d'idéal sur un anneau, elle sera développée par Dedekind [1877] : dans un anneau principal, tout idéal est engendré par son PGCD.

Pourquoi cette toute fin de son mémoire de 1764, à laquelle il semble lui-même ne pas attacher beaucoup d'importance – « je me contente d'indiquer cet usage » [Bézout 1767c] – a-t-elle donné lieu à ce que l'on appelle l'« Identité de Bézout » ?

Soit le système :

$$\begin{cases} P(x) = Ax^m + Bx^{m-1} + Cx^{m-2} + Dx^{m-3} + Ex^{m-4} + \dots\dots\dots + V = 0 \\ Q(x) = A'x^{m'} + B'x^{m'-1} + C'x^{m'-2} + D'x^{m'-3} + E'x^{m'-4} + \dots\dots + V' = 0 \end{cases}$$

où $A$, $A'$, $B$, $B'$, $C$, $C'$, etc., sont des nombres et non des polynômes.

Bézout a montré qu'il existe un polynôme non nul $L_1(x)$ de degré au plus $m'$-1 et un polynôme non nul $L_2(x)$ de degré au plus $m$-1, tels que $L_1.P+L_2.Q = 0$, si et seulement si, le Résultant du système est nul. Donc si le résultant de $P$ et $Q$ n'est pas nul, c'est à dire si $P$ et $Q$ n'ont pas de racine commune, cela signifie que le système $\mathbf{R}\begin{pmatrix} M \\ N \\ P \\ \dots \\ M' \\ N' \\ P' \\ \dots \end{pmatrix} = \begin{pmatrix} 0 \\ 0 \\ 0 \\ 0 \\ 0 \\ \dots \\ \dots \\ 1 \end{pmatrix}$ où $M$, $N$, $P$, &c.., $M'$, $N'$, $P'$, &c..sont des nombres, a une solution.

---

[31] Sylvester appelle la somme des produits de chaque polynôme par un polynôme à coefficients indéterminés, une fonction « syzygétique » des polynômes de départ [Sylvester 1853].

Par ailleurs la réciproque est évidente. On sait que dans $K[X]$, avoir une racine commune et avoir un diviseur commun sont deux propriétés équivalentes, contrairement au cas de plusieurs variables. Il est clair que, deux polynômes $P(x)$ et $Q(x)$ étant donnés, s'il existe deux polynômes $L_1(x)$ et $L_2(x)$, tels que $L_1.P+L_2.Q = 1$, alors $P$ et $Q$ n'ont pas de racine commune et pas de diviseur commun.

D'où le résultat de Bézout[32] bien connu : Les polynômes $P(x)$ et $Q(x)$ sont premiers entre eux, si et seulement si, il existe un couple unique de polynômes, $L_1(x)$ et $L_2(x)$ tels que l'on ait l'identité $L_1.P+L_2.Q = 1$, avec le degré de $L_1$ strictement inférieur au degré de $Q$, et le degré de $L_2$ strictement inférieur au degré de $P$. Si Étienne Bézout ne l'a jamais clairement énoncée, ce résultat provient assez facilement de ses méthodes de calcul de la résultante (aussi bien en 1764 qu'en 1779) et fortement suggéré dans la conclusion de son mémoire de 1764.

L'attribution du nom de Bézout à ce théorème n'a eu vraisemblablement lieu que très tard. Pendant longtemps ce résultat ne fut pas reconnu comme spécialement important. Ainsi dans les cours de la fin du XIX[e] siècle, soit il n'apparaît pas, soit on le trouve mais démontré grâce à l'algorithme d'Euclide et non attribué à Bézout[33] (voir [Alfonsi 2005, p. 347-348]). C'est grâce à la nouvelle approche avec les idéaux où la notion de PGCD rejoint celle de Bézout, que le lien a été fait par certains mathématiciens et que le théorème a pu prendre son nom. On trouve une attribution à Étienne Bézout au début du XX[e] siècle [Papelier 1903]. On le retrouve cité en 1937 [Garnier 1937], mais c'est surtout à partir de Bourbaki [1952], qui donne le nom de Bézout au résultat plus général : « Pour que les éléments $(x_i)_{i \in I}$ d'un anneau principal $A$, soient étrangers dans leur ensemble, il faut et il suffit qu'il existe des éléments $(a_i)_{i \in I}$ de $A$, nuls sauf un nombre fini d'entre eux et tels que $\sum_{i \in I} a_i x_i = 1$ », que l'identité lui est définitivement attribuée.

---

[32] Le théorème analogue dans $\mathbf{Z}$ (Deux éléments de $\mathbf{Z}$, $a$ et $b$, sont premiers entre eux si et seulement si , il existe un couple $(u,v)$ d'éléments de $\mathbf{Z}$ tels que $au+bv=1$), est souvent appelé aussi, de façon impropre, théorème de Bézout, alors qu'il est attribué à Bachet de Méziriac [Peiffer et Dahan 1986].

[33] Ni d'ailleurs à quiconque.

### IV.     La partie Algèbre du cours de mathématiques d'Étienne Bézout

Choiseul, secrétaire d'État à la Marine et à la Guerre, confie à Bézout, le 1<sup>er</sup> octobre 1764, la lourde responsabilité de réorganiser complètement les études des officiers de la marine, en le nommant examinateur des écoles de ces officiers et en lui confiant la rédaction du cours qui constituera l'essentiel de leurs études. En 1768, il lui confie les mêmes charges pour l'artillerie [Alfonsi 2005, chap.IV]. De 1764 à 1773, Étienne Bézout va donc se consacrer à la rédaction de ses cours et à sa fonction d'organisateur et d'examinateur des écoles d'officiers, fonction qui l'éloigne de Paris presque six mois par an. On comprend mieux son obligation de faire des choix dans ses sujets de recherches (voir *supra*).

Étienne Bézout a écrit deux manuels de mathématiques :

Le *Cours de Mathématiques à l'usage des Gardes du Pavillon et de la Marine* [1764-1769], qui met en lumière, par delà le cadre des écoles militaires, la réflexion très poussée de Bézout sur l'enseignement des mathématiques[34] - réflexion influencée par d'Alembert et l'*Encyclopédie*-, la haute idée qu'il se fait de l'enseignement et l'importance, pour lui, du lien entre l'enseignement et la recherche,

Le *Cours de Mathématiques à l'usage du Corps royal de l'artillerie* [1769-1772], qui est, pour la partie théorique, une version très simplifiée du premier, sauf en mécanique. C'est donc à la partie Algèbre du premier cours que nous allons nous intéresser beaucoup plus riche et innovante dans son contenu que les cours d'algèbre existant à son époque.

Le volume III, *ALGEBRE & application de cette science à l'Arithmétique & la Géométrie,* est publié en 1766. Comme on le comprend en lisant la préface, Bézout considère que c'est avec cet ouvrage qu'il va enfin pouvoir dépasser le niveau élémentaire auquel il s'est astreint dans les deux premières parties (Arithmétique et Géométrie). Il y a deux raisons à cela. La première est pédagogique, cette troisième partie est destinée à la dernière année d'études des Gardes de la Marine, donc aux plus âgés et aux meilleurs, puisqu'ils auront passé le barrage des examens. La deuxième raison lui est plus personnelle, l'analyse algébrique finie – c'est-à-dire la théorie des équations algébriques - est, nous l'avons vu, l'unique thème de recherche de Bézout à cette époque.

---

[34] Sur la question de la place des mathématiques dans les écoles militaires par opposition avec cette même place dans les autres institutions d'enseignement, et sur l'activité de Bézout pour développer encore plus les mathématiques et leurs aspects spécifiques dans les écoles de la Marine et de l'Artillerie, voir [Alfonsi 2005, chap. IV].

Quand Bézout aborde les systèmes d'équations du premier degré à deux ou plusieurs inconnues, il donne pour tous les élèves la méthode de substitution y compris pour des quantités littérales, mais il ajoute en petits caractères(c'est-à-dire en partie non obligatoire pour tous)[35], l'idée de base de son mémoire de 1764 :

« D'ailleurs nous réduirons, par la suite, l'art de chasser les inconnues dans les équations qui passent le premier degré, à celui de les chasser dans celles du premier degré. Les méthodes que l'on a eues jusqu'ici pour éliminer ou chasser les inconnues, dans les équations qui passent le premier degré, ont toutes (si l'on en excepte seulement celles qu'ont données MM. Euler & Cramer) l'inconvénient de conduire à des équations beaucoup plus composées qu'il ne faut. Ces dernières même ne sont point à l'abri de cet inconvénient, lorsqu'on a plus de deux inconnues. Il peut donc être utile de donner ici des moyens faciles pour avoir les valeurs des inconnues dans les équations du premier degré. » [Bézout 1766, p. 95].

Dans le paragraphe « Des équations à deux inconnues, lorsqu'elles passent le premier degré », Bézout innove complètement par rapport à ses prédécesseurs. Après avoir, pour tous les élèves, expliqué la méthode de substitution pour arriver à la résultante, il réserve aux meilleurs une des méthodes d'élimination de son mémoire de 1764, celle qui conduit au Bézoutien. Mais, à partir du § 166, page 205, Bézout expose sa méthode originale pour arriver à la résultante quand les deux équations ont le même degré $m$ en $x$, en y apportant une modification : il ne calcule pas le déterminant du système homogène obtenu, mais chaque variable avec les $m$-1 premières équations. Il les remplace ensuite par les valeurs obtenues, dans la dernière, qui lui donne donc la résultante cherchée. Bézout montre bien que le degré de cette résultante ne peut dépasser le produit des degrés des deux équations, ce qu'il ne faisait pas dans le mémoire.

---

[35] Cette idée - séparer nettement les connaissances de base des autres, ne rendre obligatoires que les premières mais développer le goût du savoir par les secondes - appliquée à un cours entier, constitue une innovation. Elle est matérialisée par des distinctions typographiques. Bézout peut être considéré comme le prédécesseur des auteurs actuels de cours pour les classes préparatoires, dans lesquels la même notion est souvent étudiée à deux niveaux, distingués justement par la grosseur des caractères, le niveau « Maths-sup » et le niveau « Maths-spé ».

De plus, il apporte une autre amélioration essentielle à la méthode exposée en 1764 dans le cas où les équations n'ont pas le même degré. En effet, au lieu de ce qui est proposé dans le mémoire et qui, nous l'avons vu, conduit à la méthode des coefficients indéterminés et à des calculs incertains et pénibles, il propose, pour le système suivant :

$$\begin{cases} Ax^m + Bx^{m-1} + Cx^{m-2} + Dx^{m-3} + Ex^{m-4} + \ldots\ldots\ldots V = 0 \\ A'x^n + B'x^{n-1} + C'x^{n-2} + D'x^{n-3} + E'x^{n-4} + \ldots\ldots V' = 0 \end{cases}, \qquad \text{où } n < m$$

de multiplier la deuxième équation par $x^{m-n}$ et d'appliquer alors sa méthode du Bézoutien jusqu'à ce que le multiplicateur soit de degré $n$-1, ce qui donnera $n$ équations chacune de degré $m$-1. Il continue en substituant à $x^n$, dans les diverses équations du système ainsi obtenu, sa valeur dans la seconde équation et itère cette opération jusqu'à ce que la plus haute puissance restante soit $x^{n-1}$. On a alors $n$ équations chacune de degré $n$-1, et on résout comme plus haut [*Ibid.* p. 205] en considérant le système des $n$-1 premières équations, aux inconnues $x^{n-1}$, $x^{n-2}$, $x^{n-3}$,…,$x$. Il n'emploie pas, comme dans le mémoire, les coefficients indéterminés et ramène le cas de deux degrés différents à celui de deux équations de même degré auxquelles il peut alors appliquer son procédé.

Il faut s'arrêter sur ce qui précède car il est l'exemple d'une caractéristique essentielle du Cours de Bézout : la réciprocité du lien enseignement-recherche. D'une part, les recherches récentes (mémoire de 1764) d'Étienne Bézout apparaissent dans son cours *Algèbre* de 1766, et elles l'élèvent, par-là même, bien au-dessus d'un cours classique. Il donne d'ailleurs les références de plusieurs ouvrages d'un très haut niveau, voire du niveau de la recherche de l'époque : outre les siens [Bézout 1767c, 1768], il cite aussi « l'*Analyse démontrée* du Père Reyneau », « un mémoire de Tschirnaüs publié dans les *Actes de Leipzig* », un mémoire d'Euler publié « dans le tome IX des *Nouveaux commentaires de Pétersbourg*, qui vient de paroître » [Bézout 1766, p. viij], l'*Analyse des lignes courbes* de Cramer [1750] et un mémoire de l'Académie de Berlin de d'Alembert [1748]. D'autre part, la volonté d'inscrire ces travaux dans son cours, a obligé l'auteur à les reprendre, et c'est sans doute cela qui a abouti, non seulement à une amélioration notable –le degré de la résultante obtenue par le Bézoutien est au plus égal au produit des degrés des deux équations – mais aussi à une méthode, celle du Bézoutien pour des degrés différents, qu'il n'avait pas obtenue en 1764.

Il est important de noter que Bézout a donné la primeur de ce résultat à un manuel d'enseignement et non à un mémoire académique[36]. Jusqu'à un certain point, il se comporte comme si le prestige des deux lieux de publication était le même[37]. Cela montre bien, que pour lui enseignement et recherche sont de valeur et d'importance égales, qu'ils sont aussi complémentaires et s'enrichissent mutuellement, le premier, par la clarté et la précision qu'il demande aux méthodes, obligeant la seconde à se remettre parfois en question.

Cette attitude face à l'enseignement a sans doute plusieurs causes. Bézout était un homme des Lumières qui considérait l'instruction comme primordiale pour le progrès de l'humanité [Alfonsi 2005. chap. IV et p. 371-375]. Mais la volonté de nommer des académiciens, donc des chercheurs, aux plus hautes fonctions de l'enseignement, est une volonté politique. Elle anticipe sur ce qui se créera en France, dans les grandes écoles d'abord, dans les universités ensuite : le lien enseignement-recherche. L'importance donnée dans les écoles militaires aux mathématiques, discipline dont dépendaient complètement les examens pour devenir officiers, met déjà en place la conception de la formation scientifique et technique que l'on retrouvera par la suite à l'École polytechnique. Sur beaucoup de ces points, la Révolution a donc continué et amplifié des réformes mises en œuvre sous la monarchie par des réformateurs éclairés.

---

[36] Bézout n'a jamais reparlé de ce résultat, même en 1779.
[37] Cela a nui d'ailleurs à sa postérité : Jacobi parle de la méthode de Bézout (celle du Bézoutien), non comme due à ce dernier, mais comme une méthode qu'il se « rappelle avoir lue dans le cours de Bézout et que beaucoup d'autres mathématiciens ont donnée » [Jacobi 1836, p. 297]

## V. La *Théorie générale des équations algébriques* (1779)

La période 1773-1779 est pour Étienne Bézout, celle de la consécration professionnelle et de la maturité de ses recherches mathématiques. Les connaissances, le sérieux et la force de caractère qu'il a montrés pour réformer la formation scientifique dans les écoles de la Marine, lui ont apporté la pleine confiance des ministres successifs de ce département. Il a maintenant une autorité incontestable et l'importance de son rôle pour ce corps, pleinement reconnue, va de plus lui apporter une confortable aisance. D'autre part les circonstances vont lui redonner du temps libre, ce qui va lui permettre d'être plus disponible pour sa famille et pour l'Académie des sciences, et surtout de se remettre à la recherche.

En effet, la disgrâce de Choiseul en décembre 1770 entraîne la fermeture de l'École royale d'artillerie de Bapaume[38] le 1er octobre 1772. Étienne Bézout perd le jour même ses fonctions d'examinateur de l'Artillerie et de professeur de physique expérimentale à Bapaume. Si la perte financière causée par la fermeture de l'école d'Artillerie, est lourde pour Étienne Bézout, en revanche, il va pouvoir gagner un temps appréciable : Outre les leçons et les examens à l'école de Bapaume, il n'a plus à assurer : depuis 1770, les leçons de physique qu'il donnait à Mézières (confiées à Monge, [Alfonsi 2005, chap. V]) et, depuis 1771, les examens à l'école des Gardes de Rochefort (supprimée cette année là, [*Ibid*. chap. IV]). Ses voyages, qui le mobilisaient environ six mois par an, vont donc, à partir de 1773, se réduire à un voyage annuel sur Brest et Toulon, qui lui prendra un peu moins de trois mois. Il récupère

---

[38] Le lien entre la disgrâce de Choiseul et la fermeture de l'école de Bapaume est étudié dans [Hahn 1964], [Alfonsi 2005, p. 236-237]

donc, à peu près un trimestre, sur les deux qu'il consacrait à ses tâches d'examinateur et de professeur. Par ailleurs, son travail de rédaction des cours pour la Marine et l'Artillerie est achevée depuis 1772. Il se retrouve donc, au début de 1773, beaucoup plus maître de son temps.

Depuis 1765 et son écrit *Sur la résolution générale des équations de tous les degrés*, Étienne Bézout n'a plus présenté de mémoire de mathématiques à l'Académie des sciences. Même si, on l'a vu, ses cours l'ont amené à améliorer certains de ses procédés algébriques et lui ont aussi fait faire quelques recherches touchant à la mécanique, l'optique et l'astronomie, il n'avait pas le temps et la disponibilité nécessaires pour mener à bien des travaux d'envergure. Maintenant que la situation a changé, l'emploi de ce temps retrouvé lui permet, dès 1774, de présenter à l'Académie des sciences un mémoire, malheureusement perdu[39], sur l'étude de la résultante d'un système d'équations. Ce mémoire, le premier depuis 1765, est la preuve de son retour à la recherche mathématique et plus précisément à la recherche sur la théorie de l'élimination.

C'est le 17 mars 1779[40] que Bézout présente son ouvrage, la *Théorie générale des équations algébriques*, à l'Académie des sciences, qui désigne d'Alembert, Duséjour et Laplace comme commissaires. Le rapport rendu le 17 avril 1779, est écrit de la main de Laplace. Les trois commissaires sont très élogieux et recommandent l'impression du traité, avec « approbation de l'Académie et privilège du Roi », suivant l'expression consacrée. Étienne Bézout vient présenter un exemplaire de son livre le mercredi 28 avril 1779 : c'est un très gros volume *in quarto*, de 469 pages, précédées elles-mêmes, d'une importante préface de 21 pages.

---

[39] Bézout présente ce mémoire à l'Académie en janvier 1774. Le titre apparaît dans les comptes rendus du comité de librairie, en date du 19 janvier 1774, avec la mention en marge : « imprimé ». Mais il est incomplètement recopié car sans doute trop long : « Où l'on détermine à quel degré doit monter l'équation finale résultante de l'élimination dans un nombre quelconque d'équations à &… » Les comptes-rendus des séances attestent aussi de son existence puisque l'on trouve au samedi 5 mars 1774, « M. Bézout a fini la lecture de l'écrit commencé le 19 janvier » [*RMAS* 1774, p. 69]. Malheureusement ce texte est introuvable. Étienne Bézout a, sans doute, refusé la publication et repris son manuscrit car, continuant à travailler sur le sujet, il a dû réserver son contenu pour le grand livre qu'il devait déjà avoir en vu, la *Théorie générale des équations algébriques,* publiée en 1779.
[40] Entre 1764 et 1779, des travaux d'autres auteurs ont paru sur l'élimination : [Euler 1766], [Lagrange *Oeuvres*], [Vandermonde 1776], mais traitant le sujet sous d'autres angles que Bézout, ils n'ont pas influencé ce dernier (voir [Alfonsi 2005, p. 251-253]).

Après son mémoire de 1764, Bézout, conscient de n'avoir pas résolu de façon satisfaisante le calcul du degré de la résultante pour plus de deux équations et deux inconnues, continue à travailler sur ce sujet et tient compte de ce qui lui semble avoir causé son échec. Il y fait allusion, dès 1766, dans son cours pour la Marine :

« Au reste, quoique ces méthodes auxquelles nous renvoyons [celles du mémoire de 1764], abaissent considérablement le degré auquel conduiroient celles qu'on a eues jusqu'ici, & autant qu'il est possible en n'éliminant qu'une inconnue à la fois ; il y a lieu de croire cependant, qu'il peut être encore diminué ; mais *probablement on n'y parviendra que quand on aura trouvé une méthode pour éliminer à la fois toutes les inconnues hors une* [souligné par nous] ce que je ne sache pas qu'on puisse encore pratiquer généralement sur d'autres équations que sur celles du premier degré. » [Bézout 1766, p. 208-209]

La préface est extrêmement riche. Il justifie ses recherches en analyse algébrique finie (voir *supra* l'introduction), et se livre notamment à l'analyse critique des défauts des méthodes de ses prédécesseurs, mais aussi et surtout à celle des insuffisances de son travail précédent sur le sujet (1764) :

« Je conçus dès-lors qu'en combinant les équations en plus grand nombre à la fois, on pouvoit espérer des résultats plus simples. Ce soupçon me conduisit à un travail qui a fait la matière d'un Mémoire parmi ceux de l'Académie des Sciences pour l'année 1764. Mais quoique par les moyens proposés dans ce mémoire on arrive en effet, toujours, à une équation finale beaucoup moins composée, que par les méthodes qu'on avait jusques-là, néanmoins on n'arrive pas à l'équation finale la plus simple ; & quoique le facteur qui complique le résultat soit bien moins élevé que par les autres procédés, il est en général d'autant plus composé, que les équations proposées le sont plus elles-mêmes. Ces difficultés n'ont pu que me faire sentir plus vivement combien l'analyse était encore imparfaite » [Bézout 1779 p. *viij*]

Bézout recense alors, de façon lucide, les défauts de ses précédentes tentatives :

- élimination des inconnues de manière successive ;

- recherche simultanée du degré de l'équation finale et de celui des polynômes multiplicateurs, car « ignorant pleinement quel devoit être le degré de l'équation finale, on ignoroit également celui qu'on devoit donner aux polynômes-multiplicateurs, & par conséquent aussi le nombre total de coefficiens qu'ils pouvoient fournir ; à plus forte raison ignoroit-on combien il y en avoit d'inutiles » [*Ibid.* p. *ix*].

Le passage suivant contient une analyse remarquable du déroulement de sa recherche. Bézout formulant les questions qu'il se pose et auxquelles il lui faut répondre pour continuer, décrit en détail la démarche de sa pensée :

« L'idée de multiplier les équations proposées, par des fonctions de toutes les inconnues qu'elles renferment, de faire une somme de tous ces produits, & de supposer, dans cette somme, que tous les termes affectés de toutes les inconnues qu'il s'agit d'éliminer, s'anéantissent ; cette idée, dis-je, s'étoit déjà présentée plusieurs fois à mon esprit, ainsi que probablement elle s'est offerte à d'autres. Mais quelles devoient être ces fonctions pour satisfaire à la question ? Elles pouvoient fournir moins, autant, ou plus de coefficiens qu'il n'est nécessaire pour l'anéantissement des termes à éliminer. Quel usage pouvoit-on faire des coefficiens surnuméraires ? Qui étoient-ils ? En quel nombre étoient-ils ? Et s'il étoit possible d'en employer un nombre moindre que celui des termes à faire disparoître (comme cela a lieu, en effet, dans plusieurs cas, ainsi qu'on le verra sur la fin de cet ouvrage), comment devoit-on se conduire pour ne pas arriver à des équations de condition ? Ces questions étoient précisément ce qui faisoit le noeud de la difficulté. [...]

En un mot, l'idée de procéder à l'élimination en multipliant les équations proposées restoit toujours une idée stérile, tant que ces questions n'auroient pas été résolues. » [*Ibid.*]

C'est, entre autres, sur ces questions, qu'il a buté dans le cas de trois équations à trois inconnues dans le mémoire de 1764. En effet nous avons vu qu'il n'arrivait pas à déterminer le degré des polynômes multiplicateurs et donc le nombre de coefficients avec lesquels il devait travailler.

Cette démarche analytique d'explication du déroulement de ses recherches, sans cacher *ses* [souligné par nous] échecs mais au contraire en les rendant constructifs pour ses démarches ultérieures, exposée dans un traité important, marque, elle aussi, une caractéristique de Bézout. Là aussi ses qualités didactiques, venant à la fois de sa formation, de ses responsabilités enseignantes et de la réflexion pédagogique de l'Encyclopédie [Encyclopédie, art. « élémens des sciences »] menée par d'Alembert qui fut un de ses proches, expliquent cette originalité dans l'écriture de son ouvrage.

## 1.   Le Livre I : le théorème de Bézout sur le degré de la résultante de *n* équations à *n* inconnues.

Nous allons maintenant reprendre les divers points de la démonstration du « théorème de Bézout », en respectant les notations originales de l'auteur, qui écrit :

- Un polynôme complet[41] de degré T à une seule inconnue $x$ : $(x)^T$

- Un polynôme complet de degré T à *n* inconnues u, x, y, z, etc. : $(u...n)^T$

Cette façon d'écrire les polynômes est personnelle à Bézout. Ce n'est pas l'écriture classique qu'il avait employée en 1764 et qui était la plus répandue. Sans doute, voulant ici traiter *n* équations à *n* inconnues, a-t-il pensé que l'écriture habituelle des coefficients en *a, a', a''*, etc. *b, b', b''*, etc. et le nombre trop important de termes (compte tenu de toutes les combinaisons possibles entre les exposants des variables) était impraticable et cela l'a amené à une recherche symbolique sur les notations des polynômes : ils sont caractérisés par leurs degrés et leurs nombres de variables, c'est cela qu'il a retenu pour leur écriture.

---

[41] Un polynôme complet est un polynôme dans lequel aucun terme ne manque, autrement dit, dans lequel aucun coefficient n'est nul. On rappelle que dans un polynôme à plusieurs variables, le degré du polynôme est celui du monôme de plus haut degré.

La démonstration met en œuvre le dénombrement des termes d'un polynôme et donc de ses coefficients, idée que Bézout est le premier a utiliser systématiquement dans ce domaine. Il se propose ainsi de déterminer le degré de l'équation finale résultante d'un nombre quelconque d'équations complètes, à pareil nombre d'inconnues.

Soient $n$ équations complètes à $n$ inconnues :
$$\begin{cases} (u...n)^t = 0 \\ (u...n)^{t_1} = 0 \\ (u...n)^{t_2} = 0 \\ ...... \\ (u...n)^{t_{n-1}} = 0 \end{cases}, \text{ avec } t \geq t_1 \geq t_2 \geq t_3 \geq ....... \geq t_{n-1}.$$

Il écrit, pour tout $i$ entier compris entre 1 et $n$-1, $\quad x^{t_1} = x^{t_i} + \sum_{i=1}^{n-1} (u...n)^{t_i}.A_i x^{t_1-t_i} \quad$ avec un choix des $A_i$ tels que $1 + \sum_{i=1}^{n-1} A_i a_i = 0$, les $a_i$ étant les coefficients des $x^{t_i}$ dans les différentes équations. On obtient alors $x^{t_1}$ en fonction des puissances de $x$ inférieures à $t_1$. Il fait de même pour les autres inconnues, $y^{t_2}, z^{t_3}$, etc., en prenant la valeur 0 pour l'exposant $t_j - t_i$ quand celui-ci est négatif, et en exprimant à chaque fois ces puissances en fonction des puissances inférieures de la même lettre, sans faire réapparaître les puissances disparues des lettres précédentes.

Il considère alors un seul polynôme complet de degré $T$, $Q=(u...n)^T$, dont les coefficients sont indéterminés et forme « l'équation-produit » $\quad Q.(u...n)^t = (u...n)^{T+t} = 0$.

Pour tenir compte des $n$-1 autres équations, il remplace dans cette équation-produit les termes $x^{t_1}, y^{t_2}, z^{t_3}, etc...$ par les valeurs obtenues précédemment grâce à ces équations.

Il fait le même travail pour Q. Il y a donc une première condition qui apparaît pour Q, il faut que $T \geq t_1 + t_2 + t_3 + \ldots\ldots + t_{n-1}$. Puis, grâce aux coefficients de Q, il se propose d'annuler dans l'équation $(u\ldots n)^{T+t} = 0$, tous les termes en *x, y, z, etc.* de façon à ne garder qu'une équation en *u* qui sera « l'équation résultante » du système.

Pour lui, un tel polynôme Q existe sûrement car « non seulement on conçoit que cela peut arriver ; mais on voit que cela doit arriver [...] puisque la question doit à la fin se réduire à une équation en *u*, il faut, qu'après ces substitutions, tous les termes affectés de *x, y, z, &c.,* puissent être détruits. » [*Ibid.* p. 29]

Il lui faut donc maintenant déterminer :

- le nombre de termes en *x, y, z,* etc., restant dans l'équation produit, après les substitutions envisagées ;

- le nombre de termes restant dans Q après les substitutions envisagées, qui serviront à éliminer les termes précédents.

Pour cela, il a l'idée originale de faire intervenir les *différences finies*, et le résultat cherché, le degré de la résultante, va apparaître sous la forme d'une différence finie.

### *Différences finies*

Soit *X* une fonction de *x*, il note *X'=X(x+k)*. Alors *X'-X*, la différence [finie] de *X* est notée ici *dX* ou *d(X)* ou $d(X)\ldots\begin{pmatrix} x \\ k \end{pmatrix}$, c'est la différence des valeurs de *X,* quand *x* varie de *k*.

Si *P* est une fonction à plusieurs variables, par exemple *x, y, z*, il note $d(P)\ldots\begin{pmatrix} x & y & z \\ k & l & m \end{pmatrix}$, la différence de *P*, x variant de *k*, *y* variant de *l* et *z* variant de *m*.

Puis $dd(X)\ldots\begin{pmatrix} x \\ k,k' \end{pmatrix} = d(d(X)\ldots\begin{pmatrix} x \\ k \end{pmatrix})\ldots\begin{pmatrix} x \\ k' \end{pmatrix} = d(X'-X)\ldots\begin{pmatrix} x \\ k' \end{pmatrix}$, désigne la différence seconde de *X, x* variant d'abord de *k* et ensuite de *k'*.

On peut, de la même façon, définir les différences $3^e$, $4^e$, $5^e$, etc.

Exemples[42] :

$$d(x^3 - 5x^2 + 3x - 6)...\binom{x}{k} = 3x^2k + 3xk^2 - 10xk + k^3 - 5k^2 + 3k$$

$$d[(x+a)(x+a+b).....(x+a+(n-1)b)]....\binom{x}{b} = nb(x+a+b)(x+a+2b)......(x+a+(n-1)b)$$

$$d[(x+a)(x+a+b).......(x+a+(n-1)b)]....\binom{x}{-b} = nb(x+a)(x+a+b)......(x+a+(n-2)b)$$

Il aborde ensuite le problème du dénombrement des termes en le décomposant en plusieurs étapes correspondant à autant de problèmes :

**Problème I** : *Déterminer la valeur du nombre $N(u...n)^T$, des termes d'un polynôme complet à n inconnues de degré T : $N(u...n)^T$*

On voit que $N(u...1)^T = T+1$

Puis à l'aide d'une nouvelle inconnue $x$, on rend homogène tous les termes du polynôme $(u...1)^T$, on obtient alors les termes de la dimension $T$ du polynôme $(u...2)^T$ dont le nombre est $T+1$.

Exemple donné par Bézout : « Si à l'aide de l'inconnue $y$, on rend homogènes de degré 3, tous les termes du polynôme $(u...2)^3$, c'est-à-dire tous les termes suivants

$$\begin{array}{cccc} u^3 & u^2x & ux^2 & x^3 \\ u^2 & ux & x^2 \\ u & x \\ 1 \end{array}$$

on aura les termes $u^3$ $u^2x$ $ux^2$ $x^3$ $u^2y$ $uxy$ $x^2y$ $uy^2$ $xy^2$ $y^3$ qui sont tous ceux qui peuvent composer la dimension 3 du polynôme $(u...3)^3$. » [Bézout 1779, p. 22-23]

De même, ceux de la dimension $T$-1 sont au nombre de $T$, etc., donc

$$N(u...2)^T = (T+1)+T+(T-1)+(T-2)+.....+2+1 = \frac{(T+1)(T+2)}{2}.$$

---

[42] Ces exemples sont donnés par Bézout et les deux derniers vont lui servir dans les dénombrements ultérieurs.

En recommençant la même opération, pour $(u...3)^T$ à partir de $(u...2)^T$ on trouve

$$N(u...3)^T = \frac{(T+1)(T+2)(T+3)}{1.2.3}$$

En général, par induction, il pose $\quad N(u...n)^T = \frac{(T+1)(T+2)............(T+n)}{1.2.3......n}$ .

Bézout ne fait pas de démonstration par récurrence : à son époque, l'induction était considérée comme une preuve suffisante quand la loi de formation paraissait régulière.

**Problème II :** *Dans un polynôme complet à un nombre quelconque d'inconnues u, x, y, z, etc., combien y-a-t-il de termes divisibles par $u^P$ ? combien, outre ceux-ci, sont divisibles par $x^Q$ ? combien outre les précédents sont divisibles par $y^R$ ? etc.*

*sachant que P+Q+R+etc. ≤ T, T étant le degré du polynôme.*

Le groupement de tous les termes divisibles par $u^P$ peut se mettre sous la forme $u^P (u...n)^K$ et donc K=T-P. Le nombre de termes divisibles par $u^P$ est $N(u...n)^{T-P}$.

Le nombre de termes divisibles par $x^Q$ est de la même façon $N(u...n)^{T-Q}$

et le nombre de termes divisibles par $x^Q$ outre ceux divisibles par $u^P$ est

$$N(u...n)^{T-Q} - N(u...n)^{T-P-Q} = d[N(u...n)^{T-Q}]...\binom{T-Q}{-P}$$

De proche en proche, il induit que le nombre de termes divisibles par la *(k+1)*-ième inconnue à la puissance *K,* outre les termes divisibles par les *k* inconnues précédentes à des puissances *P, Q, ..&c.,* données, est $d^k[N(u...n)^{T-K}]...\binom{T-K}{-P,-Q,..........,-J}$ où *J* est la puissance de la *k*-ième inconnue.

**Problème III :** *Le nombre de termes restant dans le polynôme* $(u...n)^T$ *une fois que l'on a enlevé tous les termes divisibles par les k premières inconnues aux puissances respectives A, B, C, &c..est* $d^k[N(u...n)^T]...\begin{pmatrix} T \\ -A,-B,.........,-J,-K \end{pmatrix}$ .

Il le montre par induction, après vérification jusqu'à $k=4$.

Si l'on reprend l'exemple précédent du *Problème I*[43], $(u...2)^3$, le nombre total de termes est $N(u...2)^3 = 4.5/2$ donc 10, comme on peut le vérifier. Le nombre de termes divisibles par $u^2$ est $N(u...2)^{3-2} = 2.3/2$, d'où 3, en effet ce sont les termes $u^3$, $u^2x$, et $u^2$. Le nombre de termes divisibles par $x$, après la suppression des termes divisibles par $u^2$, est $N(u...2)^{3-1} - N(u...2)^{3-2-1} = 3.4/2 - 1 = 5$, ce qui se vérifie puisque ces termes sont : $ux^2$, $x^3$, $ux$, $x^2$, et $x$. Donc le nombre de termes restants après avoir enlevé les termes divisibles par $u^2$ et $x$, est $10 - 3 - 5 = 2$, il ne reste effectivement que les termes $u$ et 1.

Il peut revenir à la recherche du degré de la résultante, puisqu'il sait maintenant que :

- le nombre de termes restant dans Q, après les substitutions envisagées, est

$$d^{n-1}[N(u...n)^T]...\begin{pmatrix} T \\ -t_1,-t_2,-t_3,.....-t_{n-1}. \end{pmatrix}$$

- et celui restant dans l'équation-produit après les mêmes substitutions, est

$$d^{n-1}[N(u...n)^{T+t}]...\begin{pmatrix} T+t \\ -t_1,-t_2,-t_3,.....-t_{n-1} \end{pmatrix}$$

Si D est le degré de l'équation résultante en $u$, D+1 est le nombre de termes de cette équation. Le nombre de termes à éliminer dans l'équation produit est donc

$$d^{n-1}[N(u...n)^{T+t}]...\begin{pmatrix} T+t \\ -t_1,-t_2,-t_3,.....-t_{n-1} \end{pmatrix} - (D+1),$$

qui est aussi le nombre des équations à résoudre.

---

[43] Ce que ne fait pas Bézout qui prend un exemple de polynôme complet de degré 6 à 3 inconnues, aux écritures très longues, voir [Bézout 1779, p. 27-28]

Le nombre des inconnues est celui des coefficients restants de Q moins 1, car on peut toujours, par exemple, prendre dans l'équation produit un terme égal à 1.

On doit donc avoir

$$d^{n-1}[N(u...n)^T]...\begin{pmatrix} T \\ -t_1,-t_2,-t_3,.....-t_{n-1}. \end{pmatrix}-1 \geq d^{n-1}[N(u...n)^{T+t}]...\begin{pmatrix} T+t \\ -t_1,-t_2,-t_3,.....-t_{n-1} \end{pmatrix}\text{-D-1}$$

et pour le plus petit D possible, l'égalité.

D'où la valeur de $\quad \mathbf{D = t.t_1.t_2.t_3.....t_{n-1}}, \quad$ car, d'après ce qui précède, et puisque

$$N(u...n)^{T+t} = \frac{(T+t+1)(T+t+2)............(T+t+n)}{1.2.3........n}, \text{ alors}$$

$$\mathbf{D} = d^n[N(u...n)^{T+t}]...\begin{pmatrix} T+t \\ -t,-t_1,-t_2,-t_3,.....-t_{n-1} \end{pmatrix} = \frac{1}{1.2.3......n}d^n[(T+t)^n]...\begin{pmatrix} T+t \\ t,t_1,t_2,t_3,...t_{n-1} \end{pmatrix}$$

et les deux derniers exemples des différences finies lui permettent de conclure.

Bézout énonce alors son théorème

« Le degré de l'équation finale résultante d'un nombre quelconque d'équations complettes renfermant un pareil nombre d'inconnues, & de degrés quelconques, est égal au produit des exposans des degrés de ces équations. » [*Ibid.* p. 32]

Cet énoncé appelle quelques remarques : Bézout suppose, implicitement, - puisqu'il tire la valeur de D de l'égalité entre le nombre d'inconnues et le nombre d'équations - non seulement que le système homogène dont le nombre d'équations est celui du nombre de termes à éliminer dans Q. $(u...n)^t$ et le nombre d'inconnues celui des coefficients restants dans Q après substitution, moins 1, est toujours compatible (ce qui est le cas puisqu'il est homogène). Il suppose de plus que l'équation finale est une vraie équation en *u* et non une condition toujours (ou inversement jamais), vérifiée et que les coefficients de l'équation résultante en *u* ne s'annulent pas pour le (ou les) plus haut(s) degré(s).

Or, si on arrive à une condition toujours vérifiée (en général 0=0), c'est que les équations de départ ont, toutes, un facteur commun et donc une infinité de racines communes.

Bézout n'a pas pensé à retirer ce cas des équations envisagées. Ce n'est qu'à la fin de l'ouvrage où il est confronté à un tel exemple, qu'il indique :

« On se tromperoit cependant si de ce dernier résultat on concluoit que l'une des deux équations proposées exprime toute la question » [Bézout 1779, p. 438]. Il pose, sur ce cas particulier, la question de la validité de ce qu'il a fait, et il répond :

« Non, sans doute, [ce n'est pas valide] si avant d'appliquer ce qui a été dit, on n'a pas eu soin de simplifier les équations proposées autant qu'il est possible ; c'est à dire, de leur ôter leur commun diviseur » [*Ibid.* p. 439].

Si par contre la condition n'est jamais vérifiée ni dans les réels ni dans les imaginaires, la seule solution sera infinie. Si les coefficients des plus hauts degrés s'annulent, l'équation aura autant de solutions réelles ou imaginaires que sa plus haute dimension, les autres valeurs étant infinies.

Bézout n'a pas évoqué ce dernier point mais il ne contredit pas sa conclusion.

En dehors du cas où les équations ont un facteur commun, le résultat d'Étienne Bézout est correct : les équations ne peuvent pas avoir plus de $tt_1t_2t_3......t_{n-1}$ solutions communes réelles ou imaginaires[44].

Nous avons parlé jusque là de l'énoncé.

Des limites apparaissent aussi en ce qui concerne la démonstration. En particulier, on doit constater que la « réciproque » c'est à dire la démonstration qu'il n'y a pas de racines étrangères dans la résultante trouvée et que chaque racine de la résultante correspond bien à une solution des *n* équations, ne s'y trouve pas. Bézout pense avoir évité tous les « facteurs étrangers à la question », par le remplacement dans le polynôme multiplicateur et le polynôme-produit des termes en $x^{t_1}$, $y^{t_2}$, $z^{t_3}$, etc., qui d'après lui « doivent suffire pour y exprimer toutes les conditions de la question » ; par l'usage des seuls

---

[44] En comptant les points multiples, connus à l'époque.

« coëfficiens utiles », qui mettent à l'abri d'avoir dans la résultante certaines racines du polynôme multiplicateur, puisqu'il évite « les facteurs superflus »; et par la forme générale de son procédé qui, en évitant les éliminations successives, ne peut mener qu'à une seule équation, la résultante.

Il l'écrit dans sa préface, « l'équation finale à laquelle on seroit conduit par ce procédé [des éliminations successives], peut être différente selon la manière dont on l'aura appliqué, & cependant on sent bien qu'il ne peut y avoir qu'une seule équation finale » [Bézout 1779, p. *v*]. Pour lui, la construction qu'il donne, règle d'elle-même le problème.

Bézout applique ensuite son résultat à la géométrie :
« On sait que les surfaces des corps peuvent être exprimées par des équations à trois inconnues, donc si ces corps sont tels que leurs surfaces puissent être exprimées par trois équations algébriques, il résulte immédiatement de notre Théorème général, ce Théorème général de Géométrie :*Les surfaces de trois corps dont la nature peut être exprimée par des équations algébriques, ne peuvent jamais se rencontrer toutes les trois, en un plus grand nombre de points, qu'il n'y a d'unités dans le produit des trois exposans du degré de ces équations »* [*Ibid.* p. 32]

En fonction de ce que nous avons remarqué plus haut, ce résultat est vrai en prenant la précaution de préciser que ces corps ne doivent avoir aucune composante commune (par exemple trois plans ne doivent pas avoir une droite en commun).

On a vu qu'en dimension deux, bien qu'il y ait eu des tentatives de démonstration par Euler et Cramer, la preuve de Bézout est la première vraiment satisfaisante. En dimension 3 et plus généralement en dimension $n$ quelconque, Bézout est le premier à donner une démonstration du fait que la résultante a pour degré le produit des degrés des équations.

À la fin de ce livre premier, et surtout, à la fin de son étude des équations complètes, dès la page 34, Bézout a résolu le problème essentiel qui l'avait empêché jusque-là ([Bézout 1767c]) de trouver la résultante au-delà de deux équations à deux inconnues : il connaît maintenant le degré de la résultante, ou au moins le maximum possible de ce degré, pour tout système de $n$ équations à $n$ inconnues. Il peut donc à présent, en considérant s'il le faut toutes les équations comme complètes avec des termes éventuellement nuls, déterminer le degré des polynômes multiplicateurs de chaque équation pour former la résultante, ce qu'il ne réussissait pas à faire antérieurement[45].

## 2. Le Livre II : Calcul de la résultante de $n$ équations à $n$ inconnues

Le « Livre second » de son traité va avoir dès lors pour sujet, les méthodes pour obtenir l'équation résultante de $n$ équations à $n$ inconnues. Il revient sur une des idées qui fait son originalité et qu'il a déjà exprimée en 1764 : « La méthode par laquelle, dans le livre premier, nous sommes parvenus à déterminer le degré de l'équation finale, indique assez que l'art d'éliminer, à la fois, toutes les inconnues moins une, se réduit à la méthode d'élimination dans les équations du premier degré, à un nombre quelconque d'inconnues. » [*Ibid.* p. 168]

### *Équations linéaires et déterminants*

Il expose alors une « règle générale pour calculer, toutes à la fois, ou séparément, les valeurs des inconnues dans les équations du premier degré, soit littérales soit numériques » [*Ibid.*]. Cette règle est la suivante :

« Soient *u, x, y, z, &c.* des inconnues dont le nombre soit $n$, ainsi que celui des équations. Soient *a, b, c, d, &c.* les coëfficiens respectifs de ces inconnues dans la première équation ; *a', b', c', d', &c.* les coëfficiens des mêmes inconnues dans la seconde équation ; *a'', b'', c'', d'', &c.* les coëfficiens des mêmes inconnues dans la troisième équation; & ainsi de suite. Supposez tacitement que le terme tout connu de chaque équation soit affecté aussi d'une inconnue que je représente par *t*. Formez le produit *uxyzt&c..* de toutes ces inconnues écrites dans tel ordre que vous voudrez d'abord; mais cet ordre une fois admis, conservez le jusqu'à la fin de l'opération. Echangez successivement, chaque inconnue, contre son coëfficient dans la première équation, en observant de changer le signe à chaque échange pair : ce résultat sera, ce que j'appelle, une *première ligne*.

Echangez dans cette *première ligne*, chaque inconnue, contre son coëfficient dans la seconde équation, en observant, comme ci-devant, de changer le signe à chaque échange pair ; & vous aurez une *seconde ligne*. Echangez dans cette *seconde ligne*, chaque inconnue, contre son coëfficient dans la troisième équation, en observant de changer le signe à chaque échange pair ; & vous aurez une *troisième ligne*.

Continuez de la même manière jusqu'à la dernière équation inclusivement; & la dernière ligne que vous obtiendrez, vous donnera les valeurs des inconnues de la manière suivante : Chaque inconnue aura pour valeur une fraction dont le numérateur sera le coëfficient de cette même inconnue dans la dernière ou la $n$-ième ligne, & qui aura constamment pour dénominateur le coëfficient que l'inconnue introduite $t$ se trouvera avoir dans cette même $n$.ième ligne » [*Ibid.*].

Cette règle n'est pas démontrée par Bézout et on peut supposer qu'il l'a trouvée par induction[46]. Le passage de la « dernière ligne » correspond à la règle connue aujourd'hui, d'un déterminant constitué des coefficients et des inconnues, dont on sait qu'il donne, de la même manière que la « règle générale », la valeur des inconnues.

---

[46] Comme à son habitude, Bézout donne beaucoup d'exemples d'équations à deux ou trois inconnues, à coefficients littéraux (il retrouve ainsi, sans les nommer, les formules de Cramer) ou numériques, qui illustrent bien l'application de cette règle dans le cas général, mais aussi dans des cas particuliers où elle est susceptible de simplifications.

Exemple[47] [*Ibid.* p. 176] : Soit le système $\begin{cases} ax + by + c = 0 \\ a'x + b'y + c' = 0 \end{cases}$, la *dernière lign*e de

Bézout n'est autre que le développement de $\begin{vmatrix} a & b & c \\ a' & b' & c' \\ x & y & t \end{vmatrix}$ suivant sa dernière ligne, dans lequel les coefficients de *x, y, t*, sont bien, respectivement, les numérateurs et le dénominateur des inconnues.

Bézout remarque aussi que dans le cas d'un système homogène [« sans aucun terme

absolument connu » dit-il] de *n* équations à *n* inconnues, sa « dernière ligne », donc la *n*-ième

ligne, est le déterminant du système et que l'annulation de cette *n*-ième ligne est la condition

pour avoir d'autres solutions que la solution nulle, donc la résultante du système [« l'équation

de condition nécessaire pour que toutes ces équations puissent avoir lieu à la fois » [Bézout

1779, p. 180]]. Il revient donc ici à sa conception du mémoire de 1764, dans lequel il avait

donné, pour « un nombre *n* d'équations du premier degré qui renferment chacune un pareil

nombre d'inconnues, sans aucun terme absolument connu » [Bézout 1767c, p. 292], la règle

permettant d'obtenir « la relation que doivent avoir les coëfficiens de ces inconnues pour que

toutes ces équations aient lieu »[*Ibid.*], c'est-à-dire l'annulation du déterminant. Bézout, en

1779, donne en plus, comme l'avait fait Cramer d'une autre façon, la règle pour obtenir les solutions d'un système linéaire quelconque.

En étudiant plus spécialement les systèmes linéaires homogènes, il constate que dans un de ces systèmes de n équations à (*n*+1) inconnues, l'une des inconnues étant choisie, les autres sont proportionnelles à celle-ci ;

De plus grâce à sa méthode d'écriture de lignes, il montre que l'on peut toujours trouver

« un nombre *n*+1 de fonctions d'un nombre *n*+1 de *quantités* [souligné par nous], lesquelles fonctions soient zéro par elles-mêmes » [*Ibid.* p. 181]. Pour cela il considère un système homogène, dont les $n+1$ *quantités* sont les coefficients, avec un nombre d'équations égal au nombre d'inconnues moins 1, il rajoute une de ces équations au système, et écrit que le déterminant est nul.

---

[47] Tous les exemples présentés ici sont dans le traité de Bézout.

Exemple : Soient les quantités *a, b, c, a', b', c'*, il considère $\begin{cases} ax + by + cz = 0 \\ a'x + b'y + c'z = 0 \\ ax + by + cz = 0 \end{cases}$ ,

et il écrit la nullité du déterminant, ce qui donne une égalité toujours vérifiée.

Il utilise en fait la propriété que le déterminant est nul si deux lignes sont égales, c'est à dire que le déterminant est une forme alternée.

### *Calcul de la résultante : $1^{ère}$ méthode*

Dans les « Observations générales » qui ouvrent le Livre second, Bézout revient sur le rôle des coefficients indéterminés dans son calcul du degré de la résultante et analyse la façon dont il veut les utiliser pour déterminer la résultante elle-même :

 « *Nous concevrons qu'on multiplie chacune des équations données, par un polynôme particulier, & qu'on ajoute tous ces produits. Le résultat sera ce que nous appellerons l'Équation-Somme, laquelle deviendra l'équation finale par l'anéantissement de tous les termes affectés des inconnues qu'il s'agit d'éliminer.*[souligné par nous] Il s'agit donc actuellement 1° De fixer la forme que doit avoir chacun de ces polynomes-multiplicateurs ; 2° De déterminer le nombre des coëfficiens qui, dans chacun, ne peuvent être considérés comme utiles à l'élimination ; 3° De faire connoître s'il y a un choix à faire parmi les termes qu'on doit ou qu'on peut rejetter dans chaque polynome-multiplicateur ; 4° Si on peut se dispenser de les rejetter, quel est le meilleur emploi qu'on peut en faire. » [Bézout 1779, p. 188]

Bézout vient d'établir son programme de travail, et nous allons le suivre sur un exemple de système de 3 équations à 3 inconnues, pour montrer qu'il réussit bien à résoudre des cas pour lesquels les autres mathématiciens avant lui, et lui-même en 1764, n'avaient pas abouti.

Soit le système suivant
$$\begin{cases} ax^2 + bxy + cxz + dy^2 + eyz + fz^2 + gx + hy + kz + l = 0 \\ g'x + h'y + k'z + l' = 0 \\ g''x + h''y + k''z + l'' = 0 \end{cases}$$

Il sait que la résultante sera de degré 2, qui est le produit des degrés des trois équations.

Il multiplie donc la première équation par une constante $L$, la seconde par $G'x + H'y + K'z + L'$, et la troisième par $G''x + H''y + K''z + L''$.

Ajoutant les trois produits, il obtient l'équation-somme, dont il cherche à annuler tous les termes sauf ceux ne contenant que l'inconnue $x$.

Il obtient le système suivant de 7 équations à 9 inconnues :

$$\begin{cases} Lb + G'h' + G''h'' + H'g' + H''g'' & = 0 \\ Lc + G'k' + G''k'' + K'g' + K''g'' & = 0 \\ Ld + H'h' + H''h'' & = 0 \\ Le + H'k' + H''k'' + K'h' + K''h'' & = 0 \\ Lf + K'k' + K''k'' & = 0 \\ Lh + H'l' + H''l'' + L'h' + L''h'' & = 0 \\ Lk + K'l' + K''l'' + L'k' + L''k'' & = 0 \end{cases}$$

Le nombre de coefficients inutiles à l'élimination est un. En effet :

- dans le premier polynôme multiplicateur, on ne peut faire disparaître aucun terme à l'aide de la seconde ou de la troisième équation ;

- dans le deuxième polynôme multiplicateur on peut faire disparaître un terme grâce à la troisième équation.

En récapitulant, on a donc pour l'équation-somme un terme inutile à l'élimination, dans les coefficients indéterminés introduits. On peut donc disposer arbitrairement d'un coefficient. Pour conserver la symétrie, Bézout pose une équation arbitraire, par exemple $K'h' + K''h'' = 0$.

Il a donc 8 équations à 9 inconnues. Prenant comme paramètre une inconnue autre que $L, L', L'', G', G''$, qui apparaissent dans les termes contenant uniquement $x$, il détermine ces dernières et en les substituant dans l'équation somme, il obtient la résultante :

$$L[(La + G'g' + G''g'')x^2 + (Lg + G'l' + G''l'' + L'g' + L''g'')x + (Ll + L'l' + L''l'')] = 0$$

où    $L=(k'h'')²$ ;

$L'=-[f(h'l'')-k(k'h')]h''+[e(h'l'')-d(k'l')+k(k'h')]k''$ ;

$L''=[f(h'l'')-k(k'h')]h'-[e(h'l'')-d(k'l')+k(k'h')]k'$ ;

$G'=-[c(k'h')+f(h'g'')]h''+[e(h'g'')-d(k'g'')-b(k'h')]k''$ ;

$G''=[c(k'h')+f(h'g'')]h'-[e(h'g'')-d(k'g'')-b(k'h')]k'$.

avec la notation $(ab') = ab'-a'b$, déterminant d'ordre 2, la signification étant la même pour tous les groupements de deux lettres entre parenthèses.

Il y a un facteur constant $L = (k'h'')² = (k'h'' -k''h')²$ dont il montre que l'annulation, quand elle est possible, permet d'abaisser le degré de la résultante.

L'avantage de cette méthode, selon Bézout, est qu'« il ne peut jamais se présenter de facteur qui puisse altérer le degré de l'équation finale. Le facteur ou les facteurs qui affecteront cette équation, ne peuvent jamais être que des fonctions des coëfficiens donnés des équations proposées : & ces facteurs ont, comme nous l'avons vu, l'usage important de faire connoître les cas où l'équation est susceptible d'abaissement. » [Bézout 1779, p. 335], et il souligne : « En un mot notre première méthode envisagée analytiquement, est, ce me semble, aussi parfaite qu'il est possible. » [*Ibid.*, p. 294]

Cependant, il reconnaît que cette méthode donne lieu à de très longs calculs, ce qui l'amène à en présenter une deuxième – généralisation de sa première méthode pour deux variables dans son mémoire de 1764 – beaucoup plus rapide.

### Calcul de la résultante : $2^{ème}$ méthode

Voici cette seconde façon de procéder :

Soient $n$ équations à $n$ inconnues, il choisit l'inconnue qui sera celle de la résultante, $z$ par exemple, et il factorise dans chaque équation, de façon à ce que tous les termes contenant les autres inconnues aient pour coefficient un polynôme en $z$.

On obtient alors un système de $n$ équations à $n-1$ inconnues, et l'« équation résultante est l'équation de condition entre les coëfficiens des équations proposées » .

Bézout n'a pas besoin des valeurs des coefficients qu'il introduit pour trouver la résultante, mais seulement de leur condition d'existence pour qu'ils soient non nuls.

Pour caractériser toutes les méthodes dont nous venons de parler,

- la première méthode de 1779 consiste à travailler avec toutes les variables (et toutes les équations) sans en privilégier aucune et à multiplier chaque équation par des polynômes à coefficients indéterminés contenant eux aussi toutes les variables ;

- la première méthode de 1764 pour deux équations à deux inconnues $x$ et $y$, consiste à ordonner les deux équations par rapport aux puissances d'une inconnue choisie, par exemple $x$, l'autre inconnue $y$, ne se retrouvant que dans les coefficients polynomiaux des puissances de $x$, puis à multiplier chaque équation par des polynômes en $x$ aux coefficients indéterminés ;

- enfin, la deuxième méthode de 1779, consiste à privilégier une variable, celle choisie comme l'inconnue dans la résultante, et à ordonner les équations par rapport aux autres inconnues, l'inconnue choisie ne se retrouvant que dans les coefficients polynomiaux des puissances des autres variables. Il reste ensuite à multiplier chaque équation par des polynômes aux coefficients indéterminés ne contenant que les variables non choisies.

Il compare ainsi ces deux « procédés » du traité de 1779 :

« Par le premier, jamais cette équation n'aura un degré plus élevé qu'elle ne doit l'avoir ; mais on aura toujours un grand nombre de coëfficiens à calculer, parce qu'indépendamment de l'équation finale, l'analyse donnera aux coëfficiens de cette équation finale des facteurs qui indiqueront, ou les cas de dépression, ou des solutions particulières. Par le second procédé, le calcul pour arriver à l'équation finale, sera incomparablement plus court ; il y aura beaucoup moins de facteurs ; mais ces facteurs pourront dans quelques cas, compliquer le degré de l'équation finale. » [Bézout 1779, p. 225]

Reprenons le même exemple que précédemment en le traitant par le second procédé :

Soit le système suivant :
$$\begin{cases} ax^2 + bxy + cy^2 + dx + ey + f = 0 \\ d'x + e'y + f' = 0 \\ d''x + e''y + f'' = 0 \end{cases}$$
où les 3 inconnues de départ $x$,

$y$, $z$, ont déjà été ramenées à 2, c'est à dire que les $a$, $b$, $c$, $d$, $e$, $f$, $d'$, etc., sont des polynômes

en $z$.

Puisque Bézout sait que la résultante des équations en $x$ et $y$, sera de degré 2, il choisit

pour polynômes multiplicateurs, une constante $C$ pour la première équation, $A'x+B'y+C'$

pour la seconde, et $A''x+B''y+C''$ pour la troisième.

Il obtient donc l'équation-somme

$(Ca+A'd'+A''d'')x^2+(Cb+A'e'+A''e''+B'd'+B''d'')xy+(Cc+B'e'+B''e'')y^2+(Cd+A'f'+A''$

$f''+c'd'+C''d'')x+(Ce+B'f'+B''f'+C'e'+C''e'')y+Cf+C'f'+C''f'' =0$

qui doit être le polynôme identiquement nul.

Cela donne un système linéaire de 6 équations à 7 inconnues et comme il y a un coefficient

inutile à l'élimination (voir *supra*), Bézout introduit en plus l'équation arbitraire

$B'd'+B''d''=0$.

Il écrit alors l'annulation du déterminant, nécessaire à l'existence, pour les inconnues

$C,A',B',C',A'',B'',C''$, de solutions autres que la solution nulle.

Pour cela il considère, comme il l'explique dans sa règle de formation des lignes (« première

ligne », etc., « dernière ligne », voir *supra*) $A'A''B'B''CC'C''$, expression (qui n'est pas un

produit mais seulement un jeu d'écriture), dans laquelle il remplace chaque terme

successivement par son coefficient dans les sept équations données par :

1<sup>e</sup>) coefficient de $x^2$=0 ;      2<sup>e</sup>) coefficient de $xy$=0 ;      3<sup>e</sup>) équation arbitraire ;

4<sup>e</sup>) coefficient de $y^2$ =0 ;      5<sup>e</sup>) coefficient de $x$=0 ;      6<sup>e</sup>) coefficient de $y$ =0 ;

7<sup>e</sup>) terme constant=0 ;

en changeant le signe pour les remplacements aux emplacements pairs.

Il arrive avec sa « dernière ligne », c'est-à-dire ici la septième, à l'équation :

$$(d'e'')[c(d'f'')^2+(d'e'').(de'f')-b(e'f'').(d'f')+a(e'f'')^2]=0,$$

qui est la résultante en $z$ du système homogène des sept équations (citées au-dessus) à sept inconnues, où les groupements de deux lettres entre parenthèses représentent des déterminants d'ordre 2, $[(d'e'') = d'e''-d''e' = \begin{vmatrix} d' & d'' \\ e' & e'' \end{vmatrix}]$, et les groupements de trois lettres entre parenthèses représentent des déterminants d'ordre 3, par exemple $(de'f'') = \begin{vmatrix} d & d' & d'' \\ e & e' & e'' \\ f & f' & f'' \end{vmatrix}$.

Bézout traite à part les deux cas $(d'e'')=0$ et

$$(E) \quad c(d'f'')^2+(d'e'').(de'f') - b(e'f'').(d'f')+a(e'f'')^2 =0.$$

Il remarque que le premier cas implique une diminution du degré de la résultante.

Il est conscient que cette deuxième méthode peut introduire des facteurs superflus mais il justifie cependant son intérêt par la plus grande rapidité des calculs :

« Dans la seconde méthode […] le degré apparent de l'équation finale peut dans plusieurs cas être différent du véritable. Comme les calculs, par cette seconde méthode, sont incomparablement plus courts que dans la première, l'inconvénient de rencontrer des facteurs superflus, n'est pas assez grand pour faire renoncer aux avantages qu'elle présente. Mais il est nécessaire d'avoir des moyens de dégager l'équation finale de ces facteurs, si comme il y a grande apparence, on ne peut espérer de les éviter généralement. » [Bézout 1779 p. 336]

Il explique alors comment reconnaître les facteurs superflus :

« Dans le procédé que nous avons donné, nous avons toujours un certain nombre d'équations arbitraires à former […] Comme ces équations arbitraires peuvent toujours être choisies de plusieurs manières différentes, il est clair que les variations, dans ce choix, introduiront des variations dans le facteur superflu, par conséquent dans l'équation finale apparente : en sorte que cette dernière peut toujours être regardée comme composée de deux facteurs dont l'un qui est la véritable équation finale cherchée, ne varie pas avec les équations arbitraires, & l'autre au contraire qui est le facteur superflu, varie avec ces équations arbitraires. » [*Ibid.*]

Revenons sur l'exemple traité plus haut en donnant des noms aux polynômes qui constituent les premiers membres des équations.

$$\left\{ \begin{array}{l} f(x, y) = ax^2 + bxy + cy^2 + dx + ey + f = 0 \\ g(x, y) = d'x + e'y + f' = 0 \\ h(x, y) = d''x + e''y + f'' = 0 \end{array} \right.$$

En termes actuels, la $2^e$ méthode de Bézout consiste à dire que éliminer $x$ et $y$ dans ce système revient à écrire que 6 quelconques des 7 polynômes *f, g, h, xg, yg, xh, yh*, sont liés.

Si l'on écrit l'annulation du déterminant des coordonnées des 6 premiers par exemple, dans la base *x², xy, y², x, y,* 1, on retrouve son équation (*E*).

### 3. Réception du traité

Le premier avis sur la *Théorie générale des équations algébriques* est celui des trois académiciens commissaires d'Alembert, Duséjour et Laplace, chargés de son étude pour permettre sa publication. Le manuscrit original se trouve aux Archives de l'Académie des Sciences, dans la pochette de séance du 17 avril 1779 et il est consigné dans les Procès Verbaux de l'Académie [*RMAS* 1779, f. 89-97], avec quelques erreurs de copistes. L'original est bien signé des trois commissaires, mais de la main de Laplace, ce qui signifie, d'après les habitudes de l'Académie, que c'est surtout ce dernier qui a étudié l'ouvrage.

Le rapport de Laplace est très long, puisqu'il est constitué de 15 pages manuscrites d'une écriture très serrée, et, contrairement à ceux dont l'Académie des Sciences a l'habitude - souvent un résumé des grandes lignes des ouvrages examinés -, il contient une étude détaillée et minutieuse de l'œuvre principale de Bézout. Même si Laplace reconnaît que l'exercice est difficile (« Nous allons essayer, écrit-il, de donner à l'Académie, autant qu'il est possible de le faire sans calcul, une idée du travail de cet académicien »), il donne une idée précise du contenu du traité.

Lui non plus, comme Bézout, ne voit pas le défaut de la démonstration du théorème sur le degré de la résultante – le manque de réciproque, donc de la démonstration que toute racine de la résultante donne bien une solution du système – puisqu'il écrit : « Cette forme est donc la plus simple à laquelle on puisse arriver en faisant usage de toutes les équations proposées » [*RMAS* 1779, f. 91]. Il est très admiratif de l'originalité et de l'intelligence de la méthode : « Nous remarquerons ici que c'est principalement à cette considération [utilisation du dénombrement et des différences finies] fine et importante, que M. Bezout doit l'élégance de sa méthode et la simplicité de ses résultats » [*Ibid.*], et il conclut :

« Tels sont les principaux objets que M.Bezout a discuté dans son ouvrage ; nous n'avons pû dans ce rapport, en donner qu'une idée très imparfaite, par l'impossibilité de faire entendre sans calcul des théories difficiles à saisir, lors mesme qu'elles sont présentées avec tout le développement dont elles sont susceptibles ; mais nous ne craindrons point d'estre démentis par ceux qui liront avec attention cet ouvrage, en assurant qu'il en existe très peu d'aussi utiles au progrès de l'analyse par l'importance & la nouveauté de la matièrre, & qui soient également propre à intéresser les Géomètres par la finesse et la variété des méthodes ; nous croyons donc qu'il mérite d'estre imprimé avec l'approbation et sous le privilège de l'académie. » [*Ibid.* f. 97]

On voit donc que la toute première réception de son œuvre est extrêmement favorable. En dehors de la réaction officielle de l'Académie, c'est surtout le jugement de Laplace, qui nous intéresse. Comme nous l'avons remarqué, ce jugement très flatteur n'a rien de convenu et n'est pas simplement une politesse entre académiciens, mais le résultat d'une étude attentive. Ce sentiment est d'ailleurs confirmé par la correspondance personnelle de Laplace. En effet, dans une lettre privée antérieure à sa nomination comme commissaire pour la *Théorie* *générale des équations algébriques*, puisqu'elle est datée du 25 février 1778, Laplace écrit à Lagrange qui est à Berlin : « Il ne paraît rien de bien nouveau, en Géométrie, à Paris ; mais on imprime actuellement un Ouvrage de M. Bézout, dont l'objet est une théorie générale de l'élimination entre un nombre quelconque d'équations et d'inconnues, quel que soit le degré des équations. Je ne connais cet Ouvrage que par la lecture que l'auteur en a faite à l'Académie, et par le peu qu'il m'en a dit ; il m'a paru très bon, et d'autant plus intéressant qu'il me semble que les recherches des géomètres s'étaient jusqu'ici bornées à éliminer entre deux équations et deux inconnues. » [Lagrange, *Oeuvres* t. 14, p. 80]

Son ouvrage édité, Bézout l'envoie à Lagrange, alors directeur de la classe de Mathématiques de l'Académie des Sciences de Berlin. Lagrange lui répond, le 12 juillet 1779, une lettre qui montre qu'il a étudié soigneusement l'œuvre reçue. Si son opinion est très positive, il se permet de proposer des critiques et des conseils précis :

« Monsieur, Je vous dois des remerciements infinis, non seulement pour l'honneur que vous m'avez fait en m'envoyant votre *Théorie des équations*, mais encore pour le plaisir que la lecture de cet Ouvrage m'a causé.

J'y ai trouvé beaucoup à m'instruire, et je le mets dans le petit nombre de ceux qui sont véritablement utiles aux progrès des Sciences. J'ai surtout été frappé et enchanté de l'usage que vous faites de la méthode des différences pour déterminer le nombre des termes restants, ou la différence entre le nombre des termes de l'équation somme et le nombre des coefficients utiles de tous les polynômes multiplicateurs, et pour parvenir par ce moyen à l'expression algébrique déterminée du degré de l'équation finale. *Cette partie de votre Travail est un chef-d'œuvre d'Analyse, et suffirait seule pour rendre l'Ouvrage très intéressant pour les géomètres* [souligné par nous]; mais le prix en est encore beaucoup augmenté par les autres recherches ingénieuses et savantes qu'il renferme, parmi lesquelles je distingue principalement la règle pour l'élimination dans les équations du premier degré, […] la manière d'avoir les équations de condition les plus simples au moyen des coefficients indéterminés, […] » [*Ibid*., p. 276]

Cependant Lagrange fait remarquer à Bézout que sa méthode permet de calculer facilement les inconnues autres que celle de la résultante : « Un avantage, particulier à votre méthode d'élimination et *dont vous n'avez point parlé* [souligné par nous], consiste en ce qu'on peut avoir avec la même facilité l'expression la plus simple des valeurs des autres inconnues. En effet, si D est le degré de l'équation finale en *x*, et que, dans l'équation somme, on fasse disparaître la puissance $x^D$ en conservant à sa place le terme affecté de *y* ou de *z*, etc., on aura pour la détermination de *y* en *x*, une équation de la forme[48]    $ky + (x)^{D-1} = 0$ » [*Ibid*. p. 277]

Enfin Lagrange termine en faisant constater très clairement à Bézout – bien que de façon délicate, puisqu'il introduit sa critique après avoir donné un satisfecit pour le calcul des solutions – qu'il manque, dans son travail, une partie importante, la démonstration de la réciproque de son théorème sur la résultante, lacune que nous avons déjà évoqué (voir *supra*).

L'attribution à Étienne Bézout du théorème sur le degré de la résultante de *n* équations à *n* inconnues, s'est effectuée rapidement. Il l'avait démontré pour deux équations à deux inconnues en 1764, puis dans le cas général en 1779. Ce fut Laplace qui lui attribua officiellement ce résultat en 1795, pendant les cours qu'il donna à l'École Normale de l'an III. En effet, le 11 mars 1795 (21 ventôse de l'an III) Pierre Simon Laplace donne sa sixième leçon de mathématiques aux élèves réunis dans le grand amphithéâtre du Muséum d'histoire naturelle. Après avoir expliqué comment trouver le degré de la résultante pour deux équations à deux inconnues il énonce : « Vous trouverez cette méthode exposée dans un grand détail et appliquée à un nombre quelconque d'équations et d'inconnues dans un très bon ouvrage de Bézout qui a pour titre : *Théorie des équations*. L'auteur y démontre, par une application ingénieuse du calcul des différences finies, ce théorème général, savoir que, *si l'on a un*

---

[48] Les notations sont, tout au long de la lettre, celles de Bézout, que l'on a vues au début de la démonstration de son théorème sur le degré de la résultante.

*nombre quelconque d'équations complètes entre un pareil nombre d'inconnues, le degré de l'équation finale, résultant de l'élimination de toutes les inconnues, à l'exception d'une seule, est égal au produit des degrés de toutes ces équation »* [souligné par nous] (voir [Laplace, 1992, p. 83]).

Étienne Bézout étant mort le 27 septembre 1783, il ne connut pas la reconnaissance publique de son travail, survenue douze ans plus tard.

## VI.    Conclusion

À la fin de cette étude, qu'avons-nous appris sur les résultats réels, les problématiques, les méthodes et l'originalité du travail algébrique d'Étienne Bézout, qui nous permette de mieux comprendre le mathématicien dont des outils et des théorèmes actuels portent le nom ? En savons-nous un peu plus sur sa personnalité et la façon dont il a considéré et assumé ses responsabilités d'enseignant et de pédagogue, dans le contexte de son époque, le XVIII$^e$ siècle ?

Le point le plus novateur de Bézout sur la théorie de l'élimination est, sans nul doute, son approche consistant à tout ramener à des systèmes linéaires et, lié à cela, sa conception de la résultante en tant que condition d'annulation du déterminant d'un système du premier degré. En suivant cette idée nouvelle, il a, pour deux équations à deux inconnues, d'abord transformé la méthode d'Euler de l'*Introductio in analysin infinitorum* et donné une démonstration correcte pour le degré de la résultante, puis trouvé une méthode très originale de calcul de cette résultante, celle qui conduit à l'introduction de la matrice **B**. Cependant il n'y a pas eu de réaction notable après la publication du mémoire de 1764 et la postérité de ce résultat est essentiellement due à Sylvester. Celui-ci, qui connaissait bien l'œuvre de Bézout, a repris lui aussi le problème de l'élimination sous l'angle des systèmes linéaires. Quand il présente sa méthode « dialytic », il exprime la même idée que Bézout : « In such method accordingly the process of elimination between equations of a higher degree than the first is always reduced to a question of elimination between equations which are of the first degree only » [Sylvester 1853, p. 581]. Sylvester aboutit, en suivant ce point de vue, d'une part à la matrice **R,** dont le déterminant s'appelle, on l'a vu, le Résultant (ou le déterminant de Sylvester) [Sylvester 1840], et d'autre part à la matrice **B** dont il appelle le déterminant le Bezoutiant[49] [Sylvester 1853. p. 430-431], rendant ainsi justice à l'initiateur de la méthode, ce que n'avait pas fait Jacobi en 1836. Comme l'a écrit H.S. White, président de l'American Mathematical Society : « Yet what a commentary on the futility of the best efforts is found in the fact that both Jacobi and Minding, only 60 years later, published investigations as new whose methods and results were in effect identical with Bézout's ! At least this showed not that his work was unnecessary, but only that he was in advance of his time. » [White, 1909, p. 332].

Une autre méthode originale de Bézout est celle qui consiste à considérer, certains polynômes étant donnés, l'ensemble des sommes des produits de ces polynômes par des polynômes à coefficients indéterminés. Cela lui permet, en 1764 de considérer le PGCD de plusieurs polynômes à une inconnue comme une de ces sommes (ce qui donnera naissance à l'« identité de Bézout ») et en 1779, dans le Livre II du *Traité des équations algébriques*, de chercher dans cet ensemble la meilleure résultante de plusieurs équations, c'est-à-dire celle qui n'a pas de facteurs superflus. Sylvester lui aussi travaillera avec ces objets, en appelant ces sommes de produits, des fonctions « syzygétic » [Sylvester 1853, p. 585].

On peut remarquer l'inégalité paradoxale de traitement par la postérité des deux résultats qui précèdent. La méthode de calcul de la résultante qui a conduit au Bézoutien, figurant explicitement dans son travail, passa longtemps inaperçue et fut même attribuée à d'autres mathématiciens. Son « identité », au contraire, qu'il n'a jamais énoncée clairement mais seulement évoquée sous une forme assez lointaine de celle que l'on connaît aujourd'hui, est sans doute le plus connu des résultats qui portent son nom.

---

[49] L'œuvre de Bézout n'étant pas, bien sûr, aussi impressionnante que celle d'Euler, certains auteurs ont eu tendance, quand un sujet avait été un tant soit peu abordé par ces deux mathématiciens, à toujours considérer qu'Euler en était l'initiateur. Malgré la reconnaissance de Sylvester qui l'a attribué à Bézout, comme nous l'avons vu, le Bézoutien n'a pas échappé à ce travers, puisque certains, de façon tout à fait indue, ont voulu l'attribuer à Euler (voir [Le Vavasseur 1907, p. 80-81]). D'autres, comme Jacobi, l'ayant trouvé dans un livre de cours (voir note *supra*), ont pensé qu'il s'agissait d'un résultat anonyme.

Il faut retenir aussi l'originalité de l'utilisation par Bézout de la méthode des coefficients indéterminés. Dans ses deux mémoires sur la résolution des équations, présentés à l'Académie royale des sciences en 1762 et 1765, Bézout n'innove pas, car la méthode est déjà largement connue et il l'utilise de la façon habituelle présentée par Descartes et d'Alembert. En revanche dans ses travaux sur l'élimination, aussi bien en 1764 qu'en 1779, Bézout utilise largement les coefficients indéterminés, mais, comme on l'a vu, ce n'est alors :

- ni pour égaler une « quantité » connue à une expression contenant les coefficients indéterminés (méthode de Descartes pour le $4^e$ degré, décomposition en éléments simples), puisqu'il n'y a pas d'expression ni de degré connus pour la résultante ;

- ni pour chercher un polynôme à coefficients indéterminés mais de degré fixé, qui doit obéir à des conditions (exemple d'équation différentielle donnée par d'Alembert dans l'article « coefficients indéterminés » de l'Encyclopédie, recherche de l'équation de la normale à une courbe comme le fait Descartes [1637]). En effet, comme il ne connaît pas au départ le degré de la résultante, il ne peut supposer a priori sa forme polynomiale ;

- ni **surtout** [souligné par nous] pour calculer les coefficients, comme le font Descartes, d'Alembert et les autres mathématiciens, dans leurs diverses utilisations.

Ce n'est pas le calcul de la valeur des coefficients indéterminés qui l'intéresse mais leurs conditions d'existence [Bézout 1767c, 1779] et leur dénombrement [Bézout 1779], d'où il déduit le degré de la résultante pour plusieurs équations à plusieurs inconnues et les conditions surnuméraires qui lui permettent de calculer la résultante le plus simplement possible et sans coefficients superflus. Bézout se sert donc des coefficients indéterminés uniquement comme outils de construction et non comme des paramètres de la valeur desquels dépendent ses résultats.

Dans un autre ordre d'idée, il faut remarquer la spécificité de sa manière d'introduire ses travaux, dont la préface du *Traité des équations algébriques* est un exemple caractéristique : il rappelle ses échecs passés, en analyse les raisons, en déduit toutes les questions qui en découlent et détaille le programme de travail nécessaire à l'obtention du résultat qu'il va exposer dans le traité. Il ne se contente pas d'énoncer des résultats et des démonstrations achevées, mais il fait l'effort didactique de décrire sa démarche de chercheur. Influencé par d'Alembert et l'*Encyclopédie*, il ne fait décidément pas partie de ces « maîtres de l'art, qui par une étude longue et assidue en ont vaincu les difficultés & connu les finesses » et qui « dédaignent de revenir sur leur pas pour faciliter aux autres le chemin qu'ils ont eu tant de peine à suivre. » [*Encyclopédie*, art. « Élémens des sciences » d'Alembert]. Son rôle d'enseignant et de rédacteur de cours, explique aussi cette spécificité remarquée par Bouligand qui, après avoir cité la préface du *Traité des équations algébriques* presque in extenso, conclut : « Ces importants passages prélevés dans la préface du livre de Bézout, sont très représentatifs de la pensée du savant algébriste. Mieux que tout commentaire, ils témoignent de son habileté à situer les difficultés qui se présentent et à construire en pleine conscience de la vraie raison des choses. On devine que l'auteur de pareils travaux devait être un professeur éminent. » [Bouligand 1948, p. 122]

Enfin, nous voudrions souligner une dernière caractéristique de l'œuvre de Bézout : le lien et la réciprocité enseignement-recherche, que nous avons vus se manifester très fortement dans son cours d'Algèbre[50]. Il a montré que, si un ouvrage de cours devait être au fait de la recherche pour pousser vers elle les meilleurs, réciproquement cette dernière pouvait s'enrichir des soins mêmes apportés à la clarté de l'enseignement.

---

[50] On le trouve aussi, à un degré moindre, dans la Mécanique [Bézout 1767a - b, 1772b] et le Traité de Navigation [Bézout 1769].

Comment nous apparaît finalement Étienne Bézout ? D'abord comme un mathématicien assez solitaire[51], isolé des autres académiciens des sciences par les très lourdes responsabilités qu'il avait acceptées dans les écoles militaires à partir de 1764. Ces charges l'ont de plus obligé à consacrer ses recherches à un unique thème : l'analyse algébrique finie. Ensuite comme un homme dont la modestie - que montre entre autres[52] le fait de publier des résultats de recherche dans un livre de cours – et l'isolement, ont en partie desservi la carrière de savant et d'académicien, malgré l'originalité de ses résultats, de ses méthodes et de son style de travail.

Enfin Bézout était un homme des Lumières par la haute idée qu'il avait de l'importance de l'enseignement, son travail d'examinateur, ses préfaces et ses cours le prouvent (voir [Alfonsi 2005]). Il l'était aussi par le goût certain qu'il avait pour la recherche, non seulement mathématique, mais aussi dans beaucoup d'autres domaines, où, grâce à l'Académie des sciences et à ses voyages dans les ports, il avait l'occasion de l'appliquer. L'esprit des Lumières, qui donnait pour but à l'individu de travailler pour le bien public, à l'amélioration et à la bonne marche de la société, ressort aussi dans cette phrase écrite par Étienne Bézout dans la préface de son traité de 1779 :

« Nous nous estimerons heureux si considérant le point où nous avons pris les choses, & celui où nous les amenons, on trouve que nous avons acquitté une partie du tribut que tout homme doit à la société dans l'état où il se trouve placé. » [Bézout 1779, p. *xxj*]

Pour notre part, nous nous estimerions heureuse, si, à la fin de cet article, les travaux et la personnalité de Bézout étaient un peu mieux connus et appréciés.

---

[51] Bézout semble avoir été assez isolé dans le monde des mathématiciens, surtout après 1764 : par exemple, il ne faisait partie d'aucune Académie étrangère. Alors que Bossut, si on le prend pour point de comparaison puisque Bézout et lui furent souvent en concurrence, était membre associé de la Royal Society et des Académies de Berlin et de Saint-Pétersbourg.

[52] D'autres éléments la prouvent aussi, voir [Alfonsi 2005]

# BIBLIOGRAPHIE